\documentclass[3p]{elsarticle}

\makeatletter
\def\ps@pprintTitle{%
 \let\@oddhead\@empty
 \let\@evenhead\@empty
 \def\@oddfoot{\centerline{\thepage}}%
 \let\@evenfoot\@oddfoot}
\makeatother

\usepackage{gensymb} 
\usepackage{amsmath,amssymb,amsfonts, afterpage,float,stmaryrd,bm,bbm,paralist}
\usepackage{graphicx}

\DeclareUnicodeCharacter{2212}{-}

\usepackage{subfigure}  
\usepackage{caption}
\usepackage{bbm} 

\usepackage{comment} 

\usepackage{booktabs}

\usepackage{multirow} 

\usepackage{xcolor}

\usepackage{titlesec}
\titleformat*{\section}{\large\bfseries}

\usepackage{hyperref}

\newcommand{\br}{\mathbf{r}}            
\newcommand{\bx}{\mathbf{x}}            
\newcommand{\by}{\mathbf{y}}            
\newcommand{\ve}{\varepsilon}           
\newcommand{\bigO}{\mathcal{O}}           
\newcommand{\R}{{\mathbb R}}

\newcommand{\grad}{{\nabla}}
\newcommand{\kr}{r}                       
\newcommand{\Hv}{{\mathcal H}}
\newcommand{\Ind}{\mathbf{1}}            

\begin{document}

\begin{frontmatter}
\title{Binary Level Set method for Variational Implicit Solvation Model}

\author[add1]{Zirui Zhang}
\ead{zzirui@ucsd.edu}
\author[add1]{Li-Tien Cheng}

\address[add1]{Department of Mathematics, University of California, San Diego, 9500 Gilman Drive, La Jolla, CA 92093-0112, USA}

\begin{abstract}
In this article, we apply the binary level set method to the Variational Implicit Solvent Model (VISM), which is a theoretical and computational tool to study biomolecular systems with complex topology. Central in VISM is an effective free energy of all possible interfaces separating solutes (e.g., proteins) from solvent (e.g., water). Previously, such a functional is minimized numerically by the level set method to determine the stable equilibrium conformations and solvation free energies. We vastly improve the speed by applying the binary level set method, in which the interface is approximated by a binary level set function that only takes value $\pm 1$ on the solute/solvent region, leading to a discrete formulation of VISM energy. The surface area is approximated by convolution of an indicator function with a compactly supported kernel. The VISM energy can be minimized by iteratively ``flipping'' the binary level set function in a steepest descent fashion. Numerical experiments are performed to demonstrate the accuracy and performance of our method.
\end{abstract}

\end{frontmatter}

\section{Introduction}
Water is the ubiquitous solvent of life and plays an important role in biomolecular processes such as binding and folding of protein \cite{levy_water_2006,baron_molecular_2013}. However, it remains a challenge to understand the microscopic behavior of water. A useful approach to describe the behavior of water is the explicit solvent molecular dynamics (MD) simulation, in which the configuration of individual water molecules are sampled. However, MD simulations of large biomolecular systems remains computational expensive \cite{ricci_tailoring_2018}. 

\begin{figure}[bpth]
\centering
\includegraphics[width=0.6\textwidth]{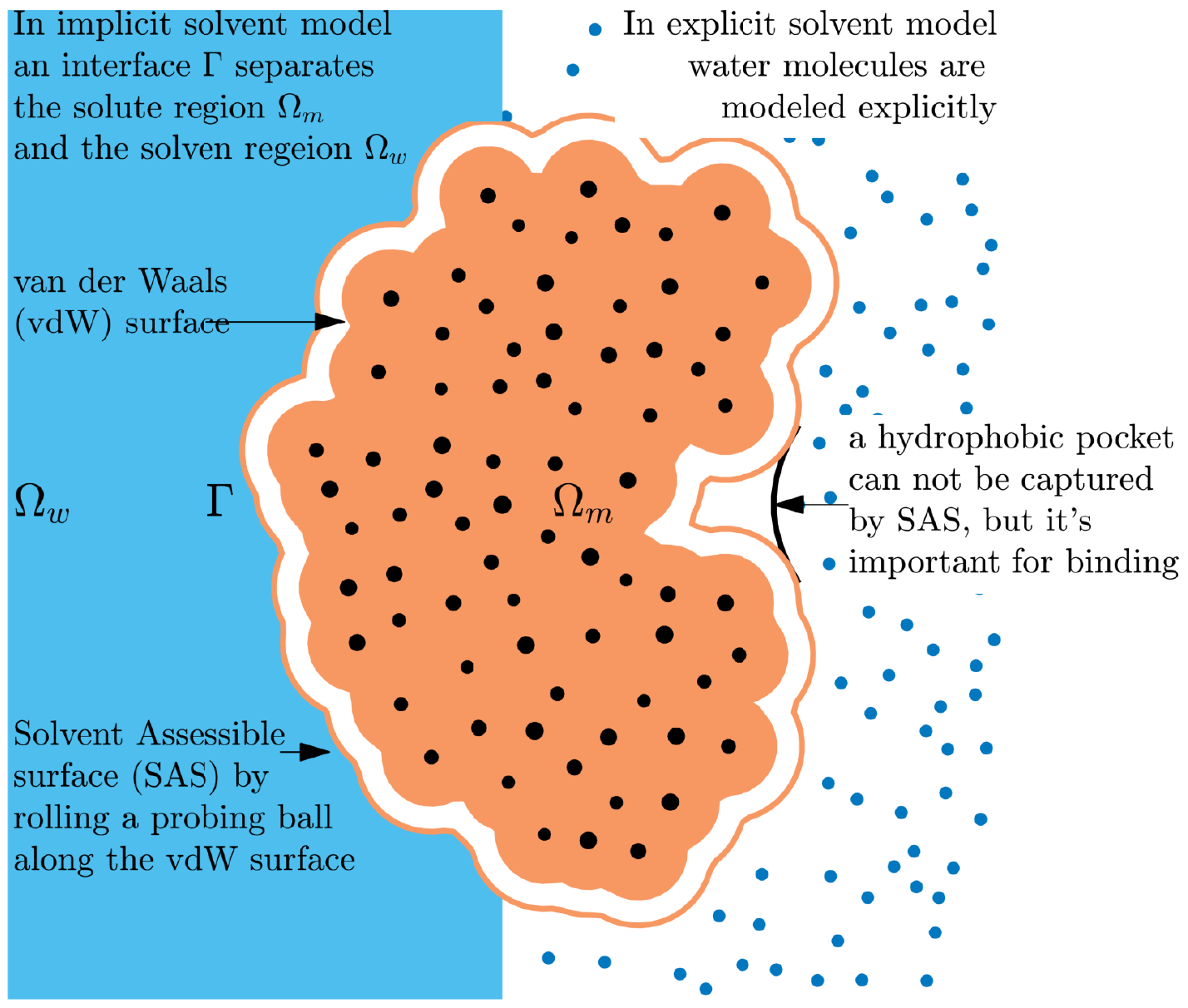}
\caption{Explicit and implicit solvent model. The solute region, solvent
region, and solute-solvent interface are denoted by $\Omega_m$, $\Omega_w$, and $\Gamma$}
\label{f:solvent}
\end{figure}

An alternative approach to MD is the implicit solvation method, which models water as continuum medium separated from the solute by an interface. Therefore, extensive statistically sampling of water configurations are no longer required. Common implicit solvation methods rely on a pre-established solute-solvent interface called solvent accessible surface (SAS), which is obtained by rolling a ball over the surface of the solute \cite{cramer_implicit_1999}. While these approaches provide useful estimation of the solvation energy in many cases, they fail to capture some important behaviors which are observed in explicit solvent simulations, such as polymodal hydration and dewetting \cite{ricci_martinizing_2017,ricci_heterogeneous_2018}. Figure \ref{f:solvent} illustrates the concepts of explicit solvent, implicit solvent and SAS.

Instead of ``guessing'' the interface, in the Variational Implicit-Solvent Model (VISM) \cite{dzubiella_coupling_2006,dzubiella_coupling_2006-1}, the solute-solvent interface is the minimizer of an energy functional defined on all interfaces enclosing the solute atoms. The functional consists of surface energy, van der Walls interaction energy, and electrostatic energy. In VISM, the contributions from different types of energies are coupled and the solute-solvent interface is the output of the theory. 

The level set method can be used to capture the stable equilibrium solute-solvent interface and compute the solvation energy \cite{cheng_application_2007,cheng_interfaces_2009,setny_dewetting-controlled_2009}. The level set method is a computational tool for tracking the motion of an interface \cite{osher_level_2003}. In the level set method, the interface is represented as the zero level set of a continuous function. The motion of the interface is governed by a partial differential equation (PDE), called level-set equation. One advantage of the level set method is that it can easily handle topological changes, such as merging and break-up. Some applications include two-phase fluid flow, crystal growth, shape optimization, image processing, etc \cite{osher_level_2003}.

The level set method (LS) had been successfully applied to numerically minimize the VISM energy functional and its various extensions.
Different formulations of electrostatic energy are incorporated, including the Coulomb Field Approximation (CFA) \cite{wang_level-set_2012} and the Poisson-Boltzmann theory \cite{zhou_variational_2014,guo_evaluation_2013}. 
One characteristic feature of VISM is its ability to capture nontrivial solvation effects such as capillary evaporation, dry-wet transition in ligand binding events, and dewetting of binding cavities \cite{zhou_variational_2014,cheng_level-set_2010, wang_level-set_2012,guo_evaluation_2013,ricci_martinizing_2017,ricci_tailoring_2018}. Stochasticity can be added to study solvent equilibrium fluctuations \cite{zhou_stochastic_2016,zhou_variational_2019}.
Solute molecular mechanics can be integrated to find the equilibrium conformation of both molecules and solvent \cite{cheng_coupling_2009}. In the LS-VISM framework, the stable equilibrium interface is obtained by evolving an initial interface in the steepest descent of the VISM energy. This approach involves solving highly nonlinear geometrical partial differential equations. Even though it's faster than the MD approach, the application of LS-VISM was limited to molecules of fixed configuration, as solving the PDE remains time consuming.

The binary level set method aims to push the limit of accuracy and performance trade-offs. In the binary level set method, the level set function only takes two discrete values, i.e. 1 and -1 \cite{lie_binary_2006,song_fast_2002,gibou_fast_2005,cheng_binary_2013}.
By ``flipping'' the value of the binary level set function between -1 and 1 \cite{song_fast_2002}, the interface can be changed drastically, opening the door for fast optimization. 
Early applications of the binary level set method focus on the piecewise constant Mumford-Shah functional in image segmentation. 
The functional consists of length of the segmenting curve and integral of image intensity in the two regions separated by the curve. This is similar to our VISM energy functional, which include the surface area of the solute-solvent interface and volume integrals in the solvent region. 
In \cite{song_fast_2002}, the flipping procedure is used to minimize the piecewise constant Mumford-Shah functional. 
Our treatment of the surface area is similar to the work in \cite{wang_efficient_2017}. 
In their work, the piecewise constant Mumford-Shah functional is minimized by an iterative thresholding method. 
The perimeter of a set is approximated by convolution of the non-local heat kernel with the characteristic function of the region.

In this work, we approximate the interface by convolution of a compact kernel with the characteristic function and arrive at a discrete formulation of the VISM energy functional. The energy can be minimized by iteratively flipping the sign of the binary level set function in a steepest descent fashion. 
Compared with the continuous LS-VISM, the binary level set method can be hundreds of times faster.
The significant speedup enables coupling VISM with intensive sampling method in molecular simulations, in which the energy needs to be evaluated more than millions of times. 
In \cite{zhang_coupling_2021}, we combine the Monte-Carlo (MC) method and binary LS-VISM to simulate the binding of the p53-MDM2 system. In this paper, we focus on the mathematical and numerical aspects of the binary LS-VISM.

The rest of the paper is organized as follows: In Section~\ref{s:theory}, we review the theory of VISM, LS-VISM and describe our binary LS-VISM. Hereafter we refer to binary level set method as BLS, and the classical level set method using continuous level set function as CLS, as we frequently compare these two approaches. Section~\ref{s:numeric} presents the numerical algorithm with implementation details. Section~\ref{s:numertest} shows convergence test and applications of our method to some molecular systems. Finally, in Section~\ref{s:conclusion}, we conclude and give an outlook to future studies. \ref{s:surfintegral} provides detailed analysis of the integral formulation of surface area of an interface represented by a binary level set function. \ref{s:surfnumeric} analyzes the discretization error of the integral formula of surface area. 
\ref{s:outside} presents a technique for integration on region outside of the computation box.

\section{Theory}
\label{s:theory}

\subsection{VISM functional}
\label{s:vism}

In VISM, the solute-solvent interface $\Gamma$ is the minimizer of an energy functional $G[\Gamma]$ defined for all possible surfaces $\Gamma$ that enclose a set of points $\br_1,\br_2,\dots,\br_N$. The energy functional includes the surface energy $G_{surf}[\Gamma] $, the van der Waals energy $G_{vdW} [\Gamma]$ with Lennard-Jones (LJ) potential $U$, and the electrostatic energy $G_{elec}[\Gamma]$ cf.~Fig.~\ref{f:schematic}.

\begin{equation}
\label{eq:vism}
G[\Gamma] = G_{surf}[\Gamma]  + G_{vdW}[\Gamma] + G_{elec}[\Gamma]
\end{equation}

The surface energy is the surface area of $\Gamma$ multiplied by a constant surface tension $\gamma_0$.
\begin{equation}
\label{eq:surf}
G_{surf} [\Gamma] =  \gamma_0 \text{Area}(\Gamma).
\end{equation}
The second term is the van der Waals (vdW) type interaction energy between the solute and the solvent
\begin{equation}
\label{eq:vdw}
G_{vdW} [\Gamma] =  \rho_w \sum_{i=1}^N\int_{\Omega_w} U_i(|\bx-\br_i|)dV.
\end{equation}
Here $\rho_w$ is the water density. $U_i$ is the Lennard-Jones (LJ) potential for atom $i$:
\begin{equation}
\label{eq:lj}
U_i(r) = 4\epsilon_i\left[ \left(\frac{\sigma_i}{r}\right)^{12} - \left(\frac{\sigma_i}{r}\right)^6 \right],
\end{equation}
where $\epsilon_i$ and $\sigma_i$ are the energy and length parameter of atom $i$.

For the electrostatic energy, we use the Coulomb Field Approximation (CFA) \cite{wang_level-set_2012}:
\begin{equation}
\label{eq:elec}
G_{elec}[\Gamma] = \frac{1}{32\pi^2\epsilon_0}\left(\frac{1}{\epsilon_w} - \frac{1}{\epsilon_m}\right)\int_{\Omega_w} \left\vert\sum_{i=1}^N\frac{Q_i(\bx-\br_i)}{|\bx-\br_i|^3}\right\vert^2dV.
\end{equation}
Here the solute atoms carry partial charges $Q_i$'s.  $\ve_0$ is the vacuum permittivity. The relative dielectric permittivities of the solute and solvent regions are denoted by $\ve_{m}$ and $\ve_{w}$. Note that $G_{vdW}[\Gamma]$ and $G_{elec}[\Gamma]$ are volume integral of the solvent region, and the integrand only depends on the positions of the atoms.

In our previous work using continuous level set method, a more complicated model for surface energy was considered, in which the surface tension is corrected by the local mean curvature of the interface. The Poisson-Boltzmann (PB) theory of electrostatics is a more accurate model and can also be incorporated into the VISM framework \cite{zhou_variational_2014}. In this work we are considering a slightly simpler model in pursuit of computational speed.

\subsection{Review of LS-VISM}
\label{s:lsvism}

In the level set method \cite{osher_level_2003,osher_fronts_1988}, the evolving interface $\Gamma = \Gamma(t)$ is represented as the zero level set of a level set function $\phi = \phi(\bx,t)$. The level set function is usually assumed to be Lipschitz continuous, and sometimes chosen to be the signed distance function to the interface. To emphasis the difference with our binary level set method, we call this classical approach the continuous level set method (CLS). In this framework, the motion of the interface is described by the PDE called the level-set equation 
\begin{equation}
   \phi_t + v_n |\grad \phi| = 0,
   \label{eq:lseqn}
\end{equation}
where $v_n$ is the velocity of the interface along the normal direction. In our previous work on CLS-VISM \cite{guo_evaluation_2013}, the interface is evolved in the direction of steepest descent of the free energy:
\begin{equation}
\begin{aligned}
   v_n(\bx) &= -\delta_\Gamma G[\Gamma]\\
   &= -2\gamma_0H(\bx) + \rho_w \sum_{i=1}^N  U_i(|\bx-\br_i|) + \frac{1}{32\pi^2\epsilon_0}\left(\frac{1}{\epsilon_w} - \frac{1}{\epsilon_m}\right) \left\vert\sum_{i=1}^N\frac{Q_i(\bx-\br_i)}{|\bx-\br_i|^3}\right\vert^2,
   \label{eq:vn}
\end{aligned}
\end{equation}
which is the negative variation with respect to local change of the interface along the normal direction. $H(\bx)$ is the mean curvature of the interface. The CLS-VISM is accurate, and it allows for more complicated VISM energy functional, such as curvature correction of surface tension \cite{wang_level-set_2012}, or Poisson-Boltzmann (PB) theory of electrostatics \cite{zhou_variational_2014}. 

Despite of its versatility and accuracy, the CLS is computational expensive because solving the nonlinear PDE \eqref{eq:lseqn} is time consuming. We need to evolve the interface from its initial shape until the steady state while we are not allowed to take large time step, because the term with mean curvature imposes stringent Courant-Friedrichs-Levy (CFL) condition of the form $\Delta t = \bigO(\Delta x^2)$, where $\Delta x$ is the grid size \cite{gibou_fast_2005,osher_level_2003}. In addition, the normal velocity depends on the number of atoms. Therefore, when going from test cases with a few atoms to practical problems involving thousands of atoms, the slowdown is significant.

\subsection{Binary level set method}
\label{s:bls}

\begin{figure}[bpth]
\centering
\includegraphics[width=0.5\textwidth]{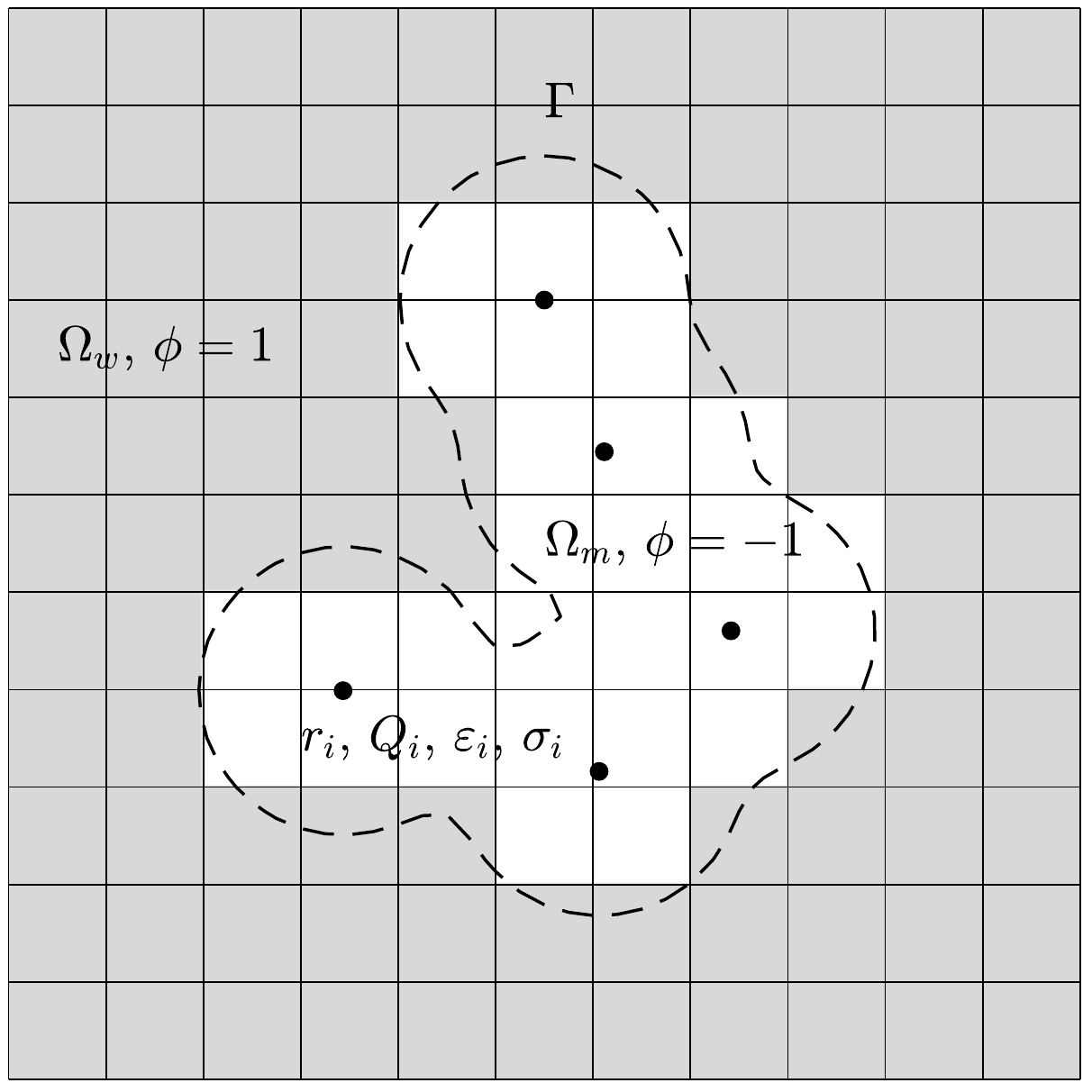}
\caption{Schematic view of a molecular system with implicit solvent. The atoms are located at $\br_i$ with charge $Q_i$ and LJ parameter $\sigma_i$ and $\varepsilon_i$. An interface $\Gamma$ (dashed line) separates the solvent region $\Omega_w$ from the solute region $\Omega_m$. In continuous level set method, $\Gamma$ is represented by the zero level set of a function. In the binary level set method, the computational domain is discretized in to grid cell. $\phi=-1$ for grid cells inside the solute region (white) and $\phi=1$ for grid cells in the solvent region (grey).}
\label{f:schematic}
\end{figure}

In the binary level set method, the interface can be approximately represented by a binary level set function $\phi$ that is $-1$ in the solute region $\Omega_m$ and $+1$ in the solvent region $\Omega_w$m, as shown in Fig.~\ref{f:schematic}. Compared with a continuous level set function, we no longer have sub-cell accuracy and the derivatives of the continuous level set function, and it's difficult to define geometric quantities such as the normal vector and the curvature. 

The early works of binary level set method focus on the piecewise Mumford-Shah functional in image segmentation. Consider a two-dimensional image $u_0$ with image domain $\Omega$. Suppose we want to approximate $u_0$ by a 2-phase image: a partition of $\Omega = \Omega_{in} \cup \Omega_{out}$ and constant intensities $c_1$ and $c_2$ in each region. The segmenting curve $\Gamma = \partial \Omega_{in}$ can be represented as the zero level set of a function $\phi$. In the binary level set framework, the following ``energy'' is minimized:
\begin{equation}
   \begin{aligned}
   F(\phi,c_1,c_2)  =  \mu \int_\Omega|\grad \Hv(\phi)|d\Omega + \lambda_1 \int_\Omega (u_0-c_1)^2\Hv(\phi)d\Omega + \lambda_2 \int_\Omega (u_0-c_2)^2(1-\Hv(\phi))d\Omega.
   \end{aligned}
   \label{eq:msf}
\end{equation}
Here, $\Hv$ is the Heaviside function. In the level set formulation, ${\rm length}(\Gamma) = \int_\Omega|\grad \Hv(\phi)|d\Omega$. $\mu$, $\lambda_1$, and $\lambda_2$ are parameters provided by the user. The 2-phase image take value $c_1$ in the region $\Omega_{in} =\{\Hv(\phi)=1\}$ and $c_2$ in the region $ \Omega_{out} = \{\Hv(\phi)=0\}$. 

We can see that the functional \eqref{eq:msf} is similar to our VISM energy functional \eqref{eq:vism}. They include the length of the curve or area of the interface respectively, and both include integral in the inside or outside region. However, in the image segmentation problem, the length serves as a regularizing term that controls the ``elasticity'' of the interface and provide a length scale for grouping object \cite{gibou_fast_2005}. For an image without noise, we can even set $\mu=0$, as the length scale is not important. However, in our energy functional \eqref{eq:vism}, the surface area needs to be approximate accurately, as it is closely relately to the surface energy.

Different approaches were developed to minimize the piecewise constant Mumford–Shah functional \eqref{eq:msf}. In \cite{lie_binary_2006}, the problem is formulated as a constrained optimization method with the constrain $\phi^2=1$, and the projection Lagrangian method and the augmented Lagrangian method are applied. In \cite{gibou_fast_2005}, the length term is first ignored, so that the Euler-Lagrange equation of \eqref{eq:msf} become an ordinary differential equation (ODE) and a larger time step can be taken. Then the length regularization is put back in by an anisotropic diffusion. The approach in \cite{song_fast_2002} is similar to our work. The Euler-Lagrange equation of the functional is ignored, and the energy is minimized directly by ``flipping'' the value of the binary level set function. However, the length of the curve is approximated in an ad-hoc way:
\begin{equation}
   \int_\Omega|\grad \Hv(\phi)|d\Omega \approx \sum_{i,j}\sqrt{ \left( H(\phi_{i+1,j})-H(\phi_{i,j}) \right)^2 + \left( H(\phi_{i,j+1})-H(\phi_{i,j}) \right)^2 }
\end{equation}
Essentially the formula approximates the curve by the edges or diagonal of the grid cells of length 1 or $\sqrt{2}$. And the formula do not converge to the length of the curve. Hence, it remains to devise a formula to approximate the surface area of an interface defined by a binary level set function.

In \cite{wang_efficient_2017,lu_threshold_2015}, the following expression is used to approximate the interfacial area between two regions:
\begin{equation}
\label{eq:perimeter}
   P_{\delta t}(\Omega) = \frac{1}{\sqrt{\delta t}}\int_{\Omega^c} G_{\delta t}* \Ind_{\Omega} dx,
\end{equation}
which is first convolving the heat kernel $G_{\delta t}$ with the indicator function of $\Omega$, then integrating in $\Omega^c$. The integral measures the amount of heat that escapes out of $\Omega$ in a short period of time, and it's known that as $\delta t$ converge to 0, the expression converge to the perimeter of a regular set \cite{jr_short-time_nodate}.

In this work, we consider a similar formula. Assuming $\Gamma$ is a smooth hypersurface in $\R^d$, then the surface area can be approximated by convolution of $\phi$ with a radially symmetric compact kernel $K$ of unit radius:
\begin{equation}
\label{eq:area}
\text{Area}(\Gamma) = C_{K,\kr,d} \int_{ \bx \in \Omega_{m}} \int_{\by \in \Omega_{w}}
K \left( \frac{|\bx-\by|}{\kr}  \right) 
\, d\by d\bx + \bigO(\kr^2) \quad \mbox{for } 0 < \kr \ll 1,
\end{equation}
where
 \begin{equation}
  C_{K,\kr,d} = \left( r^{d+1} C_d \int_0^1 K(r)r^d dr\right)^{-1} \quad C_d = \frac{2 \pi^{\frac{d-1}{2}}} {(d-1)\Gamma(\frac{d-1}{2})}
\end{equation}
Here, $\Gamma$ is the Gamma function (with a slight abuse of notation). We detail the derivation in \ref{s:surfintegral}.


\section{Numerical Method}
\label{s:numeric}
\subsection{Discretization}

In this section, we show the discretization and numerical minimization of the energy functional \eqref{eq:vism}. For notational simplicity, we assume that the energy functional \eqref{eq:vism} is restricted in some computational box $\Omega = \Omega_w \cup \Omega_m$. The computational box should be large enough to enclose the atoms and the possible interface $\Gamma$. Note that in order to find the minimizing interface, we can focus on $\Omega$, but in order to accurately compute the volume integral $G_{vdW}[\Gamma]$ and $G_{elec}[\Gamma]$, we also need to account for the integral outside of the computational box, see \ref{s:outside}.

\begin{figure}[!htpb]
\centering
\includegraphics[width=0.5\textwidth,keepaspectratio]{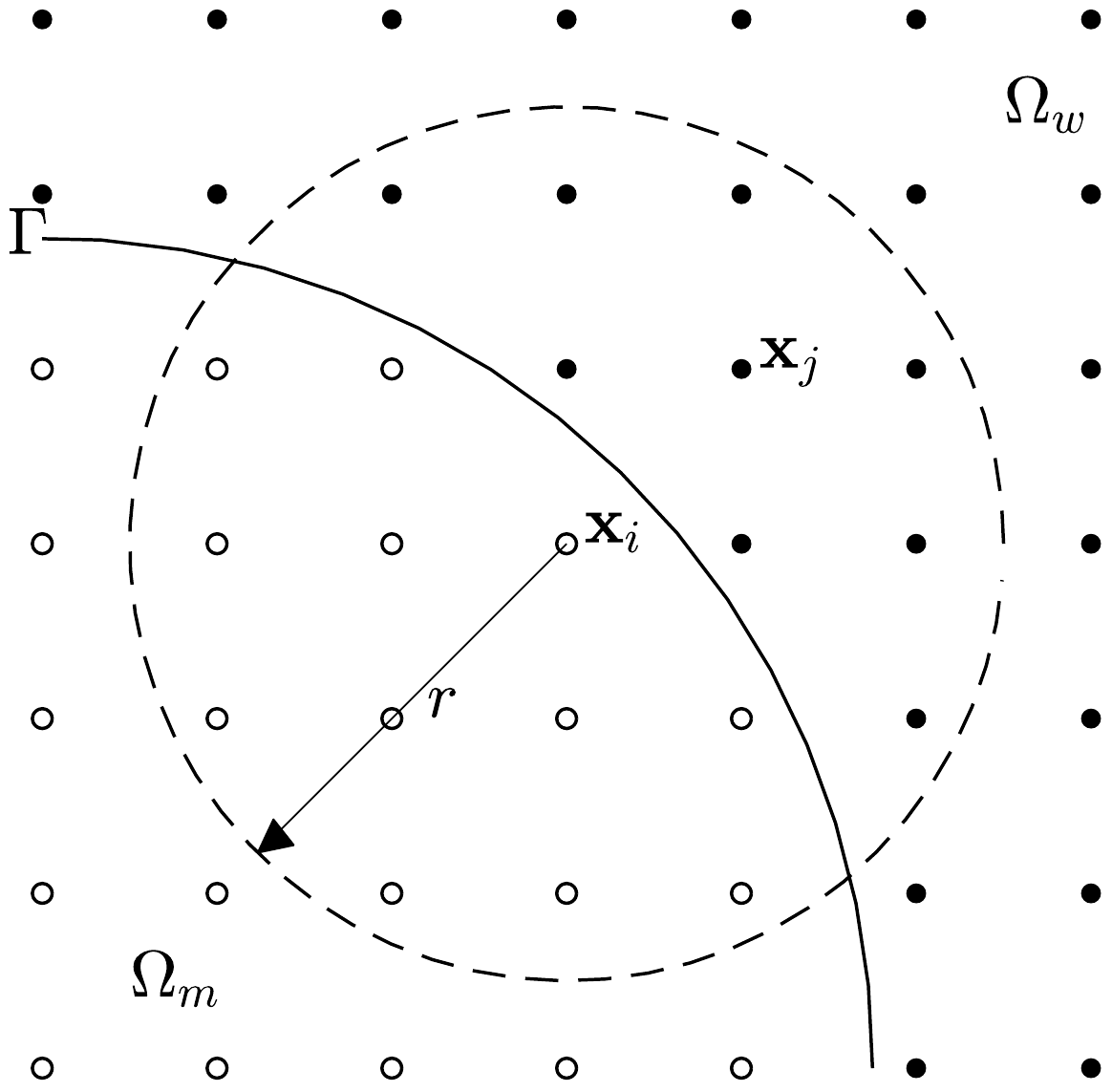}
\caption{Illustration of a scaled kernel centered at $\bx_i$ and vanishing outside a sphere (dashed line). 
Black dots represent centers of grid cells in the solvent region $\Omega_w$ and circles 
represent the centers of grid cells in the solute region $\Omega_m$.}
\label{f:kernel}
\end{figure}

With our surface area formula \eqref{eq:area}, every term in the VISM functional \eqref{eq:vism} is a volume integral and can be approximated by the midpoint rule. For the surface area \eqref{eq:area}, by choosing the kernel radius $\kr \sim \sqrt{h}$ (see \ref{s:surfnumeric}), we have
\begin{equation}
G_{surf}[\Gamma] = \sum_{\bx_i\in\Omega_m} \sum_{\substack{\bx_j\in \Omega_w\\ |\bx_j-\bx_i|\leq\kr}}K_{ij} + \bigO(h)
\label{eq:Gsurf}
\end{equation}
where 
\begin{equation}
   K_{ij} = \gamma_0 C_{K,\kr,d} h^6  K\left(\frac{|\bx_i-\bx_j|}{\kr}\right).
\end{equation}
We obtained a first order approximation of the surface area in the binary level set method. In words, we go through the grid points $\bx_i$ in $\Omega_m$, put a kernel centered at $\bx_i$, and sum up the part of the kernel in $\Omega_w$ (See Figure~\ref{f:kernel}).

The volume integrals $G_{vdW}[\Gamma]$ and $G_{elec}[\Gamma]$ can also be approximated by the midpoint rule:

\begin{equation}
G_{vdW}[\Gamma] = \sum_{x_i\in\Omega_w} G^{vdW}_i + \bigO(h)\\
\end{equation}
where 
\begin{equation}
   G^{vdW}_i = \rho_w \sum_{j=1}^{N} U_j(|\bx_i-\br_j|)h^3
   \label{eq:Gvdwi}
\end{equation}
Similarly, for the electrostatic energy
\begin{equation}
   G_{elec}[\Gamma] = \sum_{\bx_i\in\Omega_w} G^{elec}_i + \bigO(h),
\end{equation}
where 
\begin{equation}
   G^{elec}_i = \frac{1}{32\pi^2\epsilon_0}\left(\frac{1}{\epsilon_w} - \frac{1}{\epsilon_m}\right) \left\vert\sum_{j=1}^N\frac{Q_j(\bx_i-r_j)}{|\bx_i-\br_j|^3}\right\vert^2h^3
   \label{eq:Geleci}
\end{equation}
Here, we can think of $G^{vdW}_i$ and $G^{elec}_i$ as the contribution to the vdW and electrostatic energy from the grid cell which is filled with water and centered at $\bx_i$.
The discrete VISM total energy is given by
\begin{equation}
G[\Omega_m] = \sum_{x_i\in\Omega_m} \sum_{\substack{x_j\in \Omega_w\\ |x_j-x_i|<\kr}}K_{ij} + \sum_{x_i\in\Omega_w}(G^{vdW}_i + G^{elec}_i),
\label{eq:Gdiscrete}
\end{equation}
which is determined by the grid cells classified to be in the solute region $\Omega_m$.

Start with an initial guess of the surface, at each grid point, we can calculate the energy change $\Delta G_i$ if the value of $\phi$ is flipped: if $x_i\in\Omega_m$,
\begin{equation}
\Delta G_i = \sum_{\substack{x_j\in \Omega_w\\ |x_j-x_i|<\kr}}K_{ij} - \sum_{\substack{x_j\in \Omega_m\\ |x_j-x_i|<\kr}}K_{ij} + G^{vdW}_i + G^{elec}_i.
\label{eq:deltaGi}
\end{equation}
Here, the first two terms compute the difference of the kernel in $\Omega_w$ and $\Omega_m$. The third and the fourth terms are the contribution to the vdW energy and electrostatic energy if we fill the cell with water. If $\bx_i\in\Omega_w$, the change will be $-\Delta G_i$. Note that if $\bx_i$ is flipped, then the $\Delta G_j$ of the neighboring grids ($\bx_j$ with $|\bx_j - \bx_i| < \kr$) needs to be updated because of the compactness of the kernel. More specifically, suppose we flip one grid cell at each step of our algorithm, and let $\Delta G_i^{(k)}$ be the energy change if $\phi_i$ is flipped at the step $k$, then 
\begin{equation}
   \Delta G_j^{(k+1)} = \begin{cases}
   -\Delta G_j^{(k)}   & j = i \\
   \Delta G_j^{(k)} - \phi_i \phi_j 2 K_{ij}  & j\neq i, |\bx_j - \bx_i| < \kr  \\
   \Delta G_j^{(k)} & \text{otherwise}\\
   \end{cases}
\label{eq:update}
\end{equation}
Notice that the work to perform the update is proportional to the number of grid cell covered by the kernel.



\subsection{Algorithm with computational details}

\begin{figure}[!htpb]
\centering
\includegraphics[width=0.5\textwidth,keepaspectratio]{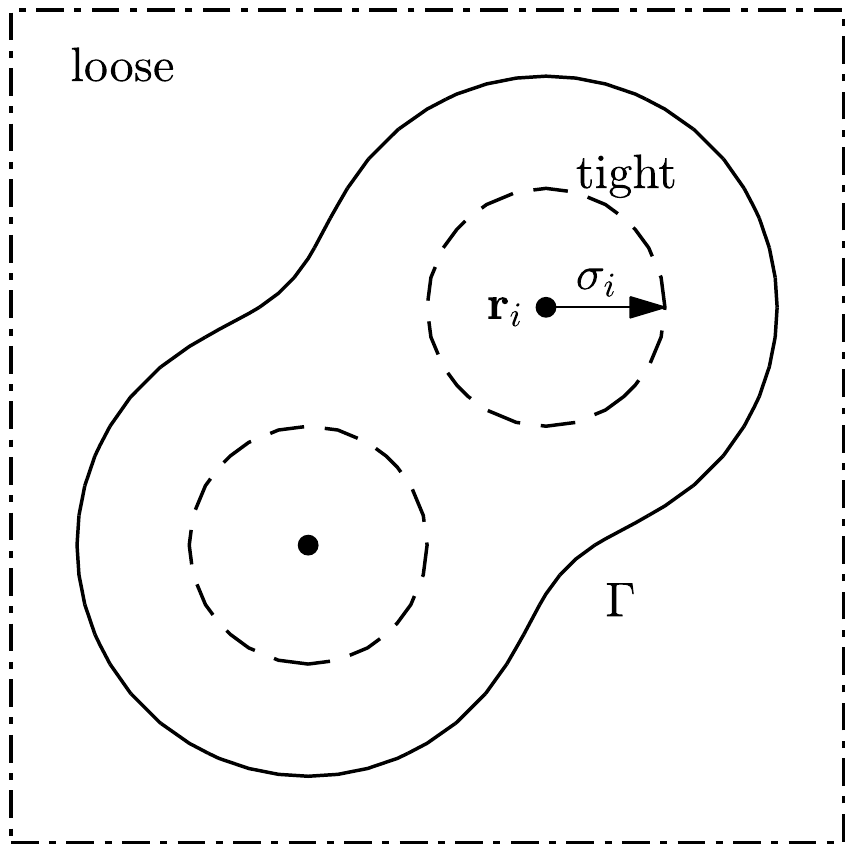}
\caption{Illustration of a tight initial surface (dashed line), a loose initial surface (dash-dotted line), and a VISM relaxed surface (solid line) $\Gamma$ surrounding the atoms (dots). For loose initial, we set $\phi=-1$ for all grid cells. The tight initial is the union of sphere of radius $\sigma_i$ centered at $\br_i$. }
\label{f:initial}
\end{figure}

The discrete VISM energy functional \eqref{eq:Gdiscrete} can be optimized in a steepest descent fashion by finding the grid cell with the largest energy decrease and flipping its $\phi$-value. Due to the non-convexity of the VISM energy functional, different initial surface may relax to different local minimizer. These local minimizers correspond to polymodal hydration states. In order to capture different local minimizer, we usually use two types of initial surface: one is a tight wrap that is the union of spheres centered at the solute atoms, the other is a large surface that loosely encloses all the atoms. In our calculation, we take the tight initial to be the vdW surface, which is the union of sphere of radius $\sigma_i$ centered at $\br_i$. For the tight initial, we mark all the grid point in the computational box as solute region (see Figure~\ref{f:initial}).

Since we needed to repeatedly find the grid cell with the largest energy decrease, we can use the min heap data structure, which takes logarithmic time to remove the smallest element and insert an element. The element of the heap is the pair ($\Delta G_i$, $i$).
We implement our binary heap with some extra information to make our flipping algorithm fast, at the expanse of extra memory use. 
We only keep elements with negative energy change and indices next to the interface, as the interface evolves during the flipping procedure. The aims is to keep the number of elements in the heap as small as possible. 
As a result, points may be inserted into the heap or deleted several times during the whole flipping procedure. 
Also notice that $\Delta G_i$ can be computed in time proportional to the number of points in the kernel, but updated in constant time as in \eqref{eq:update}. 
Therefore, as the interface evolves, we will compute and store $\Delta G_i$ whenever the $i$-th grid cell is first encountered, and update its value afterwards, even though the element may not be in the heap.

At the beginning, we insert all the interface points with negative energy change into the heap. Suppose we first remove ($\Delta G_i$, $i$) and flip the sign of $\phi_i$. Then we go through all the neighboring grid points $x_j$ within the kernel, and there are a few different situations:
\begin{compactenum}
\item If the $x_j$ is not an interface point before the flip, but a new interface point afterward, then we compute and store $\Delta G_j$. 
\item If the $x_j$ is an interface point before and after the flip, we update $\Delta G_j$ using \eqref{eq:update}. 
\item If the $x_j$ is an interface point before the flip, but no longer an interface afterward, we delete the element from the heap.
\item If the $x_j$ is not an interface point before and after the flip, nothing is performed.
\end{compactenum}
Here, $\Delta G_j$ is stored whenever it is computed or updated, but only inserted into the heap if it is negative. 

Next, we present the outline of our algorithm.
\begin{compactenum}
\item[{Step 1.}]
Input all the parameters
$\gamma_0$, $\rho_w$, $\ve_0$, $\ve_m$, $\ve_w$, 
and atomic parameters $\br_i,$ $Q_i$, $\varepsilon_i,$ and $\sigma_i,$ for all $i = 1, \dots, N$.  
Choose a computational box $[-a,a]^3$ according to the atomic coordinates and discretize 
the box uniformly with the prescribed computational grid size $h$ or the number of intervals $n$ on each side of the computational box. Note that $h = 2a/n$.
Initialize the kernel function and the binary level set function with either loose or tight initial.

\item[{Step 2.}] Compute and store $G^{vdW}_i$ (Eq.~\eqref{eq:Gvdwi})
and $G^{elec}_i$ (Eq.~\eqref{eq:Geleci}) at centers $\bx_i$ of all grid cells. 
These components do not change throughout the binary level-set iteration, and they can be parallelized.

\item[{Step 3.}] Compute and store $\Delta G_i$ (Eq.~\eqref{eq:deltaGi}) 
of grid cell near the interface. Insert the pair ($\Delta G_i$, $i$) 
to the heap data structure if $\Delta G_i<0$ and $i$ is next to the interface.

\item[{Step 4.}] Extract and remove the minimum element ($\Delta G_i$,$i$), flip $\phi_i$, and update $\Delta G_j$ at the neighboring center point $\bx_j$ with $|\bx_j - \bx_i| \leq \kr$. Modify the heap accordingly (see previous discussion).

\item[{Step 5.}] Repeat Step 4 until $ \Delta G_i > 0 $ for all grid cells. 

\end{compactenum}

As the flipping proceeds, the total energy of the system decreases monotonically. In the end, we reach a local minimum where there is no single flipping that can decrease the energy. 
However, there might be simultaneous flipping that can further reduce the energy. 

There are different ways to carry out the flipping procedure. In Jacobi iteration, one goes through all the grid cell and flips all the grid cells with negative $\Delta G$ together, and then compute their updated $\Delta G$. In a Gauss-Seidel iteration, one flips the grid cell one by one and update $\Delta G$ immediately. Here we are taking a flow-based approach, aiming to find the local minimum that is close to the initial guess. From our experience, for simple shapes, these different approaches lead to the same minimum.

\section{Numerical Experiments}
\label{s:numertest}

In this section, we present several numerical tests, all in three dimensions. We first test our formula of approximating surface area \eqref{f:area} using a sphere. We experiment with different kernel functions and kernel radii. Then we test our method with different systems: one atom, two atoms, and two protein pairs. 
We mainly compare the new BLS-VISM with the previous CLS-VISM, and demonstrate that the BLS-VISM can achieve similar accuracy as CLS-VISM but in much shorter time.

For CLS-VISM calculations, the forward Euler method is used to discretize the time derivative in the level-set equation, and a fifth order WENO (weighted essential-non-oscillation) scheme is used to discretize the spacial variables. For details, the reader can refer to \cite{wang_level-set_2012}.

For comparison of speed, we only compare the time for flipping in BLS-VISM and solving the time dependent PDE in CLS-VISM. We exclude the time for initialization, which comprises step 1 and 2 in the algorithm, and is common in both BLS-VISM and CLS-VISM. All the tests are performed on a 2017 iMac with 3.5 GHz Intel Core i5 and 16GB memory. 
We denote $n$ the number of sub-intervals in each dimension of the cubical computation box. So the total number of grid points is $(n+1)^3$. 

\subsection{Approximation of Surface Area}

\begin{figure}[htbp]
\includegraphics[width=0.5\linewidth,keepaspectratio]{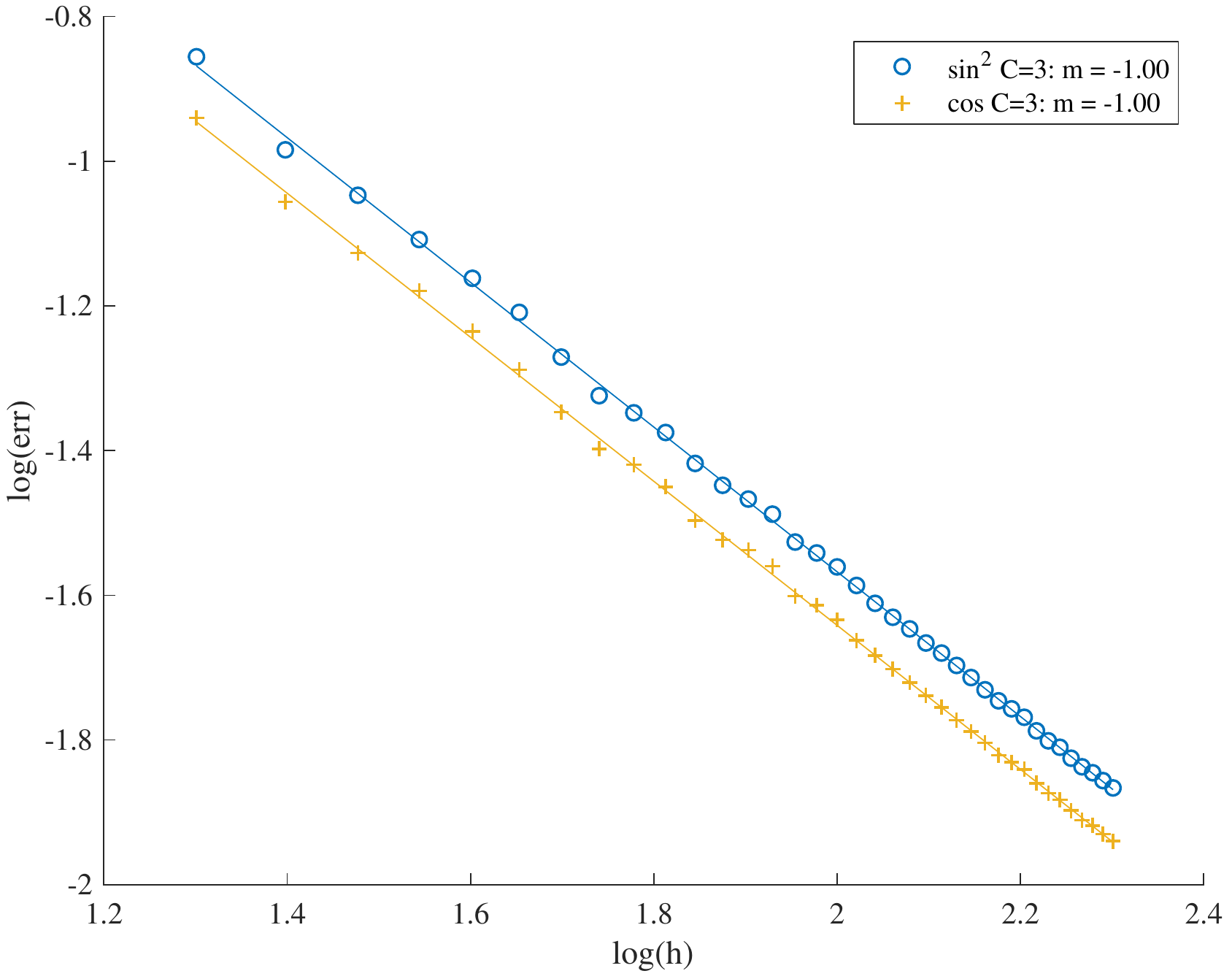}
\includegraphics[width=0.5\linewidth,keepaspectratio]{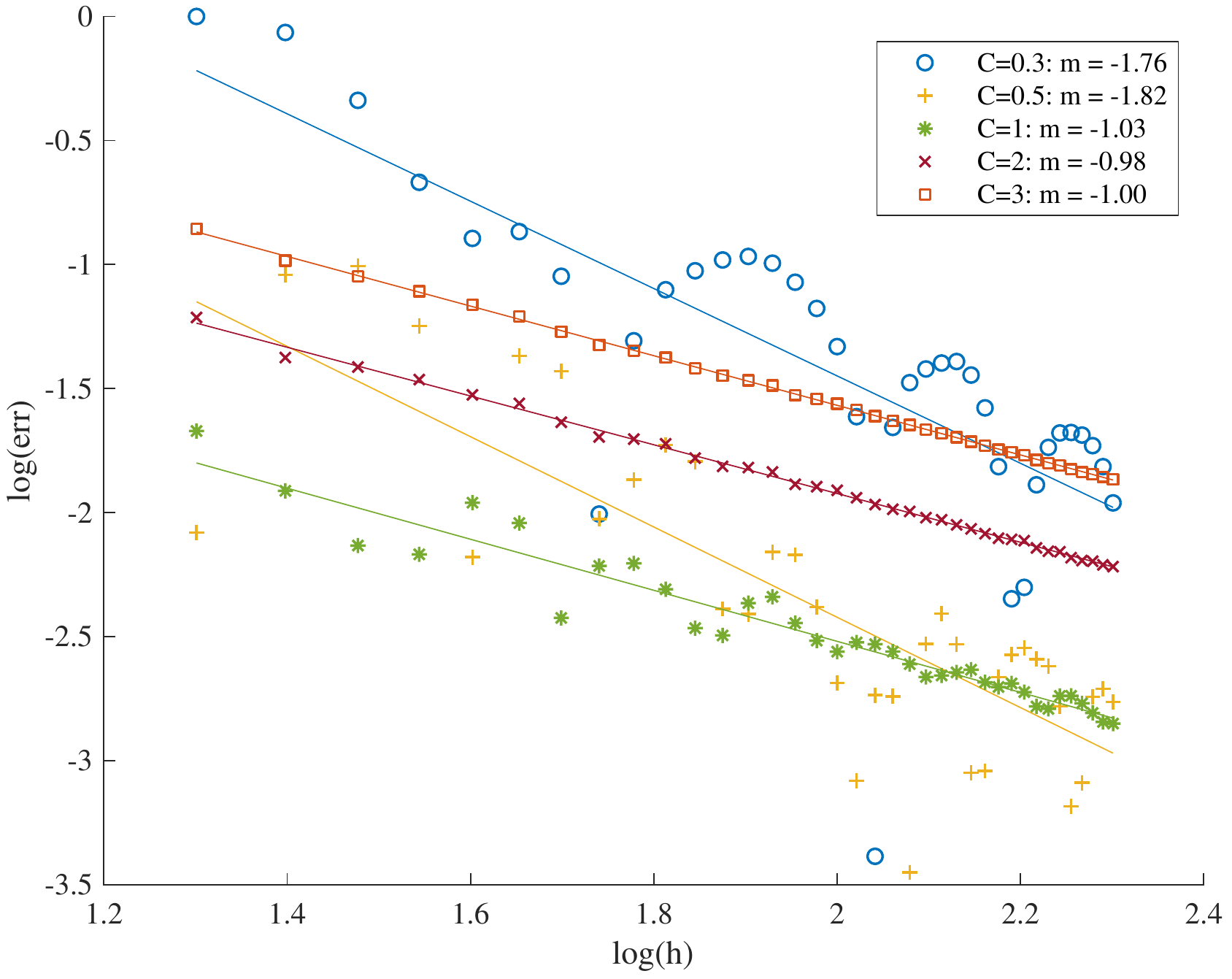}

\caption{log-log plot of the relative error versus the number of interval $n$ in one edge of the computational box. $n$ ranges from 20 to 200 with increments of 5. Each data point is an average of 6 spheres with random centers. $m$ is the slope of the line fitted by least-square. (L) sin-squared kernel and cos kernel with the same size. (R) sin-squared kernel with different size parameter $C$.}
\label{f:area}
\end{figure}

Here we apply our method to approximate the area of a sphere. 
Consider a sphere of radius 0.5 in a $[-1, 1]^3$ box. Let $\kr = C \sqrt{h}$. Here we consider 2 different kernel functions: 
$K_1(r) = \sin(\pi r)^2$ and 
$K_2(r) = \cos(\pi r)+1$ for $0\leq r\leq 1$. 
Both $K_1$ and $K_2$ are $C^\infty$ functions. The parameter $C$ control the size of the kernel. In Fig. \ref{f:area}, we briefly look at the performance of different kernel and the effect of $C$. In both figures, we plot the relative error versus the number of sub-intervals $n$ in one edge of the computational box as $n$ ranges from 20 to 200 with increments of 5. Each data point is an average of 6 trials in which the center of the sphere is randomly perturbed.

In Fig. \ref{f:area} (L), we fix $C=3$ and use different kernel functions. The approximations are all first order accuracy as the slopes of the least-square lines are close to 1. 
In Fig. \ref{f:area} (R), we use the sin-squared kernel function $K_2(r)$ and look at the effect of different values of C. The parameter $C$ controls the number of grid points in the kernel, which is proportionate to the amount of work for updating the energy change \eqref{eq:deltaGi} or computing the surface area. In three dimensions, the number of grid point in the kernel is of order $\bigO((C/\sqrt{h})^3)$. Therefore, a smaller value of C is preferable in terms of computation speed. However, as shown in Fig. \ref{f:area} (R), for $C=5,3,1,0.5,0.3$ in order, the overall error first decreases and then increases, while variance of the convergence line keep increasing. For $C=5$ or $3$, the data is fitted nicely by a straight line with slope close to 1. However, for $C=0.5$ or $0.3$, the data oscillate wildly as $n$ increases. However, the overall accuracy is still at least first order. Therefore, for demonstration of clear and stable convergence behavior, a relatively larger value of C is used.

In practice, we rarely send $n$ to infinity. Instead, we usually focus on a relatively small range of resolution, which is limited by our computation resources and the requirement on accuracy. Also the length scale is dictated by the physics of the problem at hand. For example, in our application of BLS-VISM to protein binding simulation in Section \ref{s:protein}, the average radius of the atoms in a protein is around 3 {\AA} and we only consider the grid size around 1 {\AA} for reasonable accuracy and computation time. Then we would choose the parameter $C$ based on numerical experiments with similar setting.

\subsection{One atom}

We consider a single charged atom carrying a partial charge $Q$. 
The VISM free-energy functional Eq.~\eqref{eq:vism} is then a function of the radius $R$ of that spherical 
solute region:
\begin{equation}
\label{eq:1dvism}
G(R) = 
4\pi R^2\gamma_0 +16\pi\rho_w \ve\left(\frac{\sigma^{12}}{9R^9} -\frac{\sigma^6}{3R^3}\right) 
+ \frac{Q^2}{8\pi\ve_0 R}\left(\frac{1}{\ve_w}-\frac{1}{\ve_m}\right),      
\end{equation}
where $\sigma$ and $\ve$ are the LJ parameters between the atom and a water molecule. The function $G(R)$ can be minimized accurately. We use the following parameters that are close to real systems \cite{cheng_application_2007}: $\gamma_0 = 0.174 k_BT/{\rm \AA}^2$, $\ve_m = 1$, $\ve_w=80$, $\rho_0=0.0333{\rm \AA}^3$, and $Q=1$.





In Table~\ref{t:1atom}, we show the convergence of total energy and its different components. The computation box is $[-5,5]^3$(unit {\AA}), the number of interval $n$ in one edge of the computational box goes from 20 to 200, with increment of 5. The y axis is the relative error with respect to the exact solution \eqref{eq:1dvism}. The results are averaged over 6 trials in which the position of the atom is perturbed randomly near the origin. The initial guess is a sphere of radius $\sigma$ (tight fit). The kernel radius is chosen to be $\kr = 5\sqrt{h}$. We can see that the total energy is first order accurate, so is the surface energy. The vdW energy $G_{vdW}(R)$ and the electrostatic energy $G_{elec}(R)$ exhibit slightly higher order of convergence but more oscillation. These two components of the energy are sensitive to the location of the interface, as the $G_{vdW}(R)$ increases to infinity at order $\bigO(R^{-9})$, and $G_{elec}(R)$ goes to negative infinity at order $\bigO(R^{-1})$.

\begin{figure}[htbp]
\centering
\includegraphics[width=0.7\linewidth,keepaspectratio]{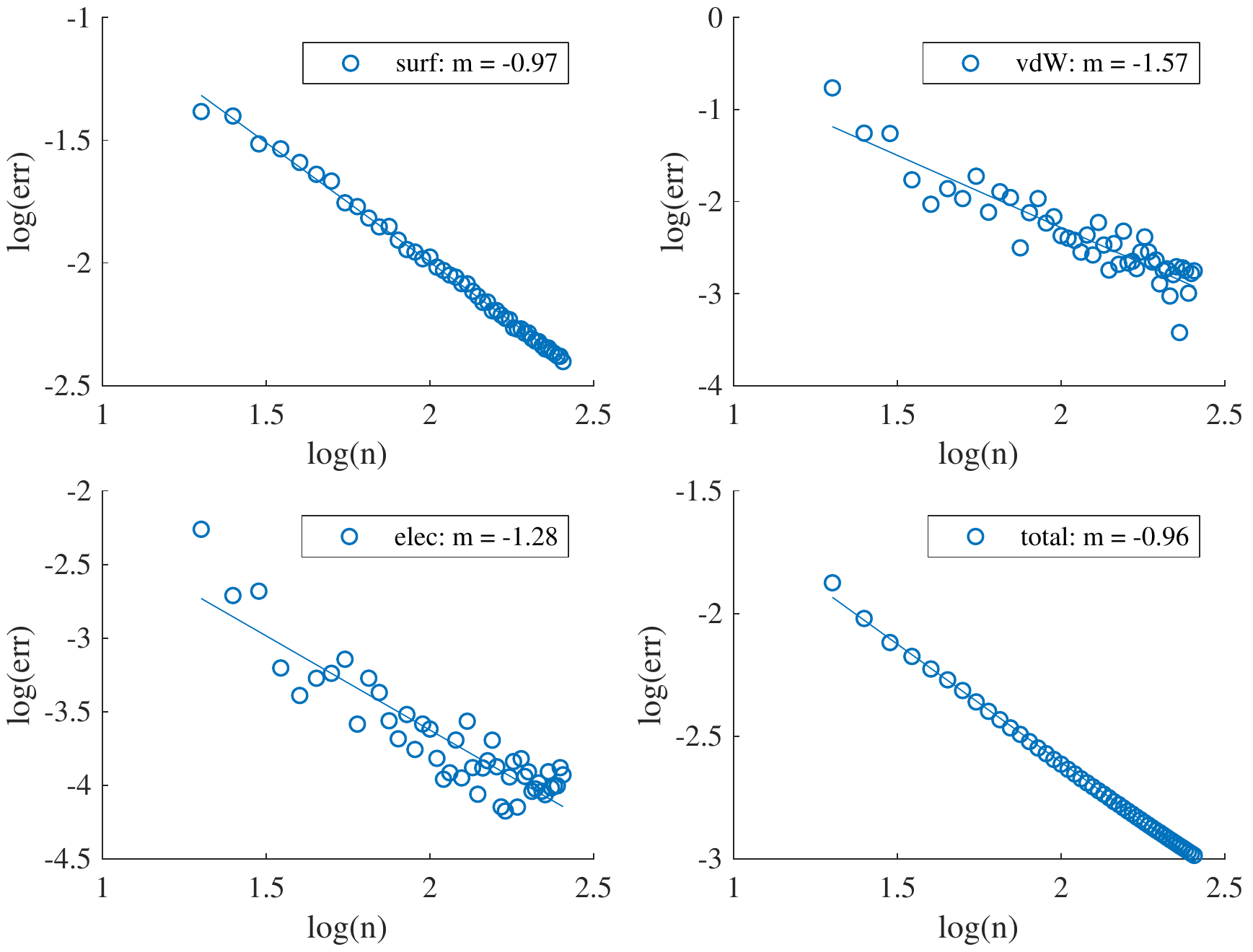}
\caption{log-log plot of the relative error of each component versus the number of interval $n$ in one edge of the computational box. $n$ ranges from 20 to 200 with increments of 5. $m$ is the slope of the line fitted by least-square}
\label{f:1atom}
\end{figure}

In Table~\ref{t:1atom}, we show a comparison of the calculation speed between the CLS-VISM  and the BLS-VISM. 
The computation box is $[-5,5]^3$(unit {\AA}) and the atom is fixed at the origin. The numbers in the parenthesis are the relative error with respect to the exact solution \eqref{eq:1dvism}. Both methods are very accurate in terms of total energy, the total error is less than 1\% with n = 50. The total energy exhibit first order convergence for BLS and second order for CLS. The individual component might have larger relative error. One reason is that the vdW and electrostatic energy are sensitive to the boundary, as mentioned before. Another reason is that the energy components have different magnitudes. In this example, the electrostatic is much smaller in magnitude. Hence even though the relative error may seem large, it has insignificant impact on the error of total energy. What's remarkable is that the BLS-VISM is about 100 times faster while still accurate for our application.



\begin{table}[!htpb]
\centering
\caption{Solvation free energy ($k_BT$) and computation time ($s$) for different grid numbers. Here, cont.\ stands for the continuous LSM and binary for the binary LSM. The number in the parenthesis is the relative error with respect to the exact solution \eqref{eq:1dvism}}
\label{t:1atom}
\medskip
\resizebox{\textwidth}{!}{
\begin{tabular}{|c|c|c|c|c|c|c|c|c|c|c|}
\hline\hline
\multirow{2}{*}{n }&\multicolumn{2}{c|}{Surf} &\multicolumn{2}{c|}{vdW} &\multicolumn{2}{c|}{Elec} &\multicolumn{2}{c|}{total} &\multicolumn{2}{c|}{Time}\\
\cline{2-3} \cline{4-5}\cline{6-7} \cline{8-9} \cline{10-11}
 & CLS & BLS & CLS & BLS & CLS& BLS& CLS& BLS& CLS& BLS\\
\hline
$25^3$ & 17.0(0.9\%) & 17.7(5.4\%) & -97.1(1.9\%) & -96.2(-2.9\%) & 2.6(48.0\%) & 3.9(-24.3\%) & -77.5(0.6\%) & -74.6(-3.3\%) & 0.006& 3.7 \\
$50^3$ & 16.8(0.1\%) & 17.1(1.8\%) & -98.3(0.7\%) & -98.2(-0.8\%) & 4.2(17.1\%) & 4.6(-9.4\%) & -77.4(0.3\%) & -76.5(-0.8\%) & 0.043& 7.1 \\
$100^3$ & 16.7(0.4\%) & 16.9(0.7\%) & -98.9(0.1\%) & -98.6(-0.4\%) & 4.9(4.0\%) & 4.8(-6.3\%) & -77.2(0.2\%) & -77.0(-0.2\%) & 0.977& 103.9 \\
$200^3$ & 16.7(0.4\%) & 16.9(0.5\%) & -99.0(0.0\%) & -98.8(-0.2\%) & 5.1(0.8\%) & 4.8(-5.9\%) & -77.2(0.1\%) & -77.1(-0.0\%) & 24.055& 2952.1 \\
\hline
\end{tabular}
}
\end{table}

\subsection{Two atoms}

In this experiment, two atoms are placed at $(-d/2,0,0)$ and $(d/2,0,0)$ for $d$ = 4, 6, 8 in the usual \textit{xyz} coordinate system. The physical parameters are the same as those in the previous experiment. 

The computation box is $[-10,10]^3$ (unit {\AA}). We choose $n=100$ for BLS and $n=50$ for CLS. Because BLS is only first order accurate, a finer grid is needed to achieve similar accuracy as the CLS. In the binary LS-VISM, we choose the kernel radius $\kr = 3\sqrt{h}$. For both continuous and binary LS-VISM, the tight and loose initial relax to the same final interface and energy. Notice that topological changes are handled easily as the center-to-center distance increases. 

\begin{figure}[!htpb]
\centering
\includegraphics[width=0.7\linewidth,keepaspectratio]{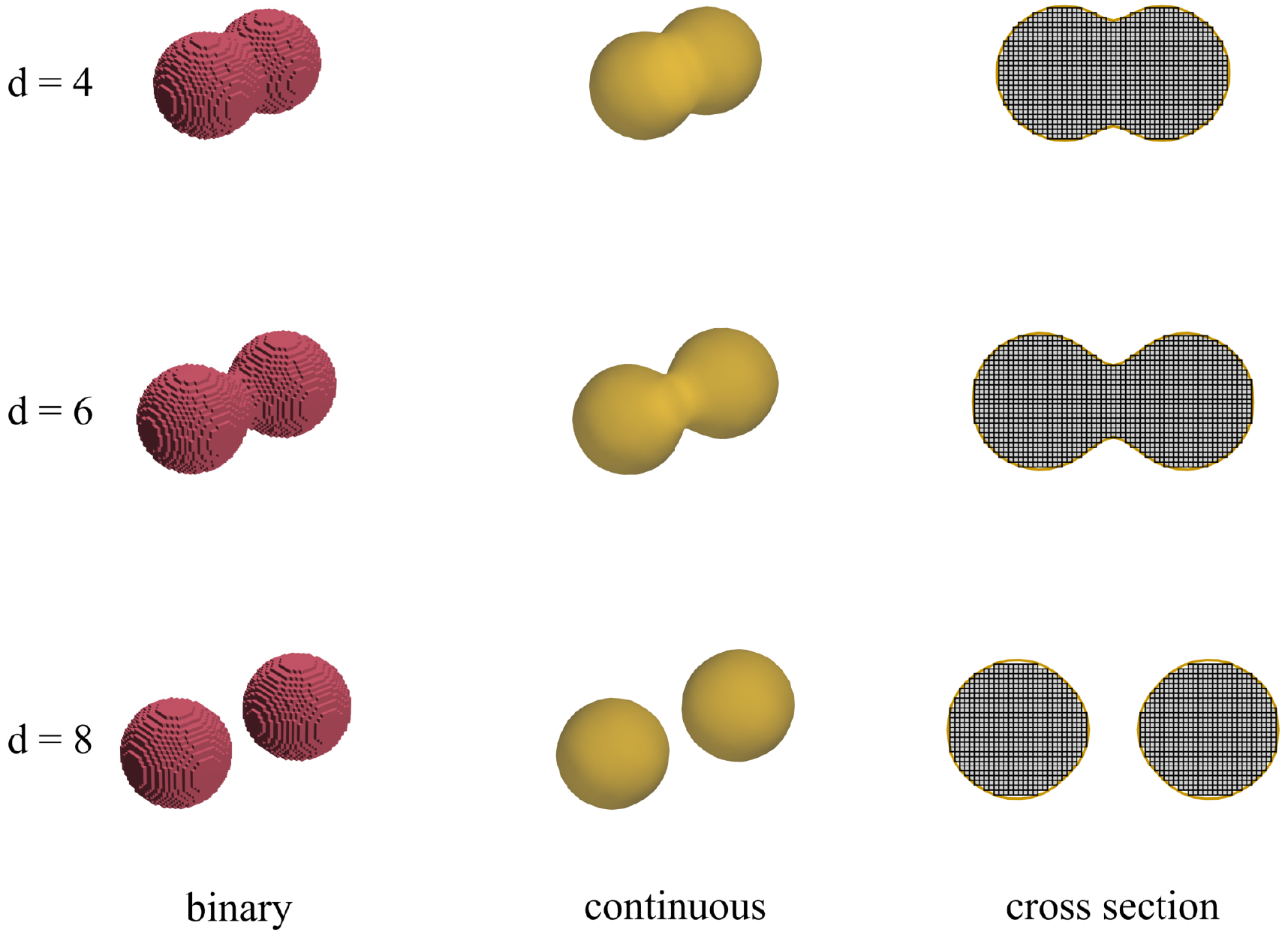}
\caption{Final interfaces of d = 4, 6, 8 \AA. Left: BLS-VISM. Middle: CLS-VISM. Right: cross section of both interface at $z=0$. $n=100$ for BLS and $n=50$ for CLS. }
\label{f:2atomsurf}
\end{figure}

In Fig.~\ref{f:2atomsurf}, we show the final interfaces from BLS and CLS for different configurations. We also show the overlaid cross sections of both interfaces. For the interfaces from BLS, we render the faces of the voxels, which are the little cubes of side length $h$ and centered at the grid points with $\phi = -1$. We can see that the interfaces from BLS closely approximate those from CLS.

In Table~\ref{t:2atom}, we show a comparison of the energy and speed between the BLS and CLS, both using tight initial. Again, BLS-VISM can obtain accurate energy and it's hundreds times faster than CLS-VISM.

\begin{table}[!htpb]
\centering
\caption{Solvation free energy ($k_BT$) and computation time ($s$) for the two atoms system. $n=100$ for BLS and $n=50$ for CLS.}
\label{t:2atom}
\medskip
\resizebox{\textwidth}{!}{
\begin{tabular}{|c|c|c|c|c|c|c|c|c|c|c|}
\hline\hline
\multirow{2}{*}{d }&\multicolumn{2}{|c|}{Surf} &\multicolumn{2}{c|}{vdW} &\multicolumn{2}{c|}{Elec} &\multicolumn{2}{c|}{total} &\multicolumn{2}{c|}{Time}\\
\cline{2-3} \cline{4-5}\cline{6-7} \cline{8-9} \cline{10-11}
&BLS & CLS & BLS & CLS& BLS& CLS& BLS& CLS& BLS &CLS\\
\hline
4& 27.0 & 27.8 & 14.1 & 12.2 & -317.9 & -314.3 & -276.8 & -274.2 & 0.7 & 11.4 \\
6& 32.2 & 33.2 & 10.5 & 9.8 & -285.9 & -283.1 & -243.3 & -240.1 & 0.7 & 20.7 \\
8& 33.3 & 34.6 & 10.1 & 8.8 & -265.6 & -262.4 & -222.2 & -219.0 & 0.8 & 5.7 \\
\hline
\end{tabular}
}
\end{table}

\subsection{Biomolecules}
\label{s:protein}

We apply our method to two complex biomolecular systems from the Protein Data Bank (PDB) \cite{berman_protein_2000}: p53-MDM2 (PDB ID 1YCR) \cite{kussie_structure_1996} and BphC (PDB ID 1DHY)\cite{senda_three-dimensional_1996}.
The p53-MDM2 system consists of more than a thousand atoms. It is a relevant pharmacological target for anticancer therapeutics \cite{chene_inhibiting_2003}.
During their binding process, the binding pocket fluctuates between dry and wet states \cite{guo_heterogeneous_2014}.
The BphC is a enzyme and contains more than four thousands atoms. The parameters for the atoms come from the force field in CHARMM36\cite{huang_charmm36m_2017}. For both experiments, to generate the positions of the atom, we start with the bounded structure from PDB, and manually pull away the protein pairs along their center-to-center axis for a distance $d$. Here $d=0$ corresponds to the bounded state.

\begin{figure}[!htpb]
\centering
\includegraphics[width=1\linewidth,keepaspectratio]{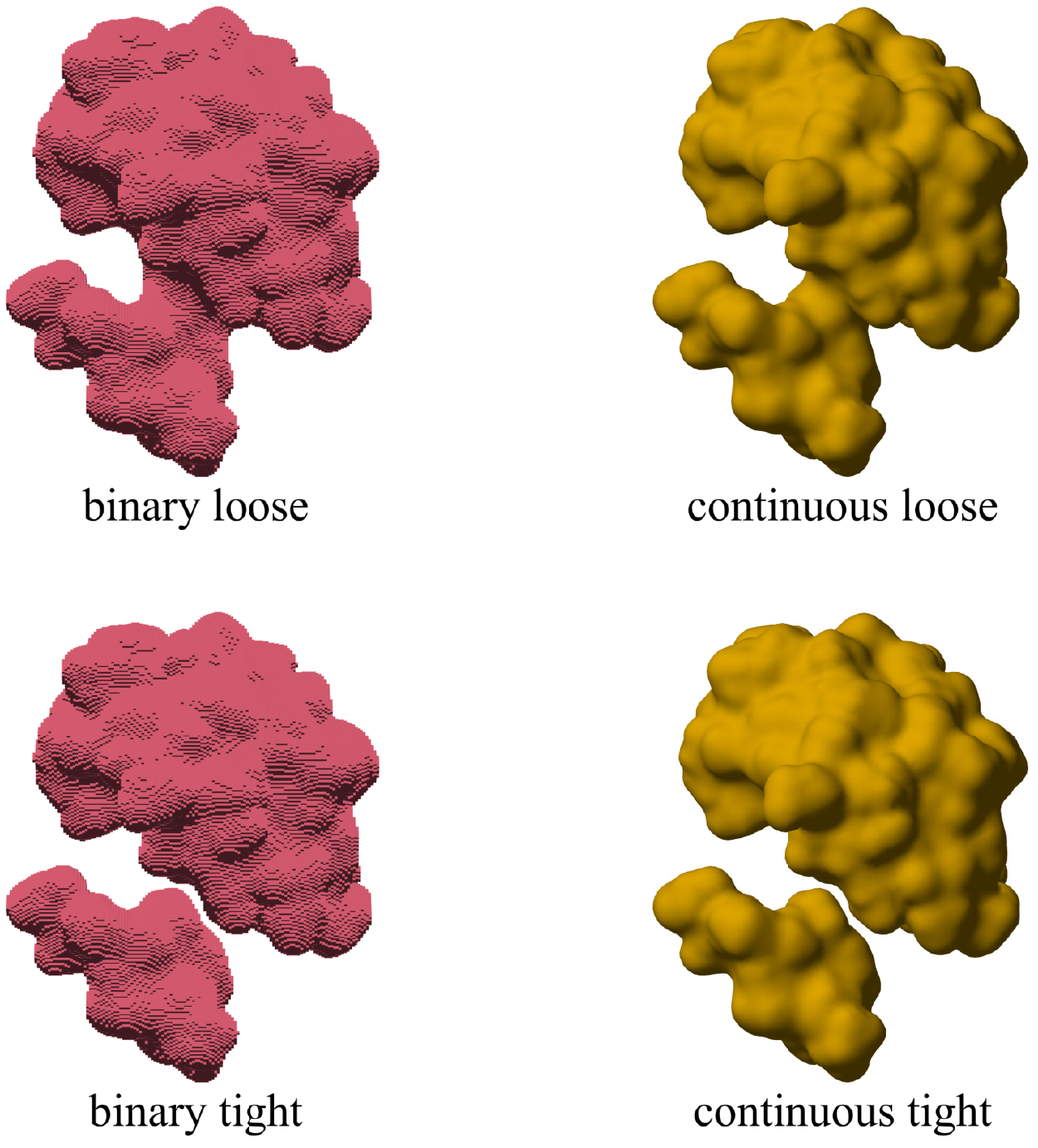}
\caption{Stable equilibrium solute−solvent interfaces of p53-MDM2 obtained at d = 14 {\AA}. The computation box is $[-31.7945,31.7945]^3$ and $n=200$ for BLS-VISM, $n=100$ for CLS-VISM}
\label{f:p53mdm2}
\end{figure}

\begin{figure}[!htpb]
\centering
\includegraphics[width=1\linewidth,keepaspectratio]{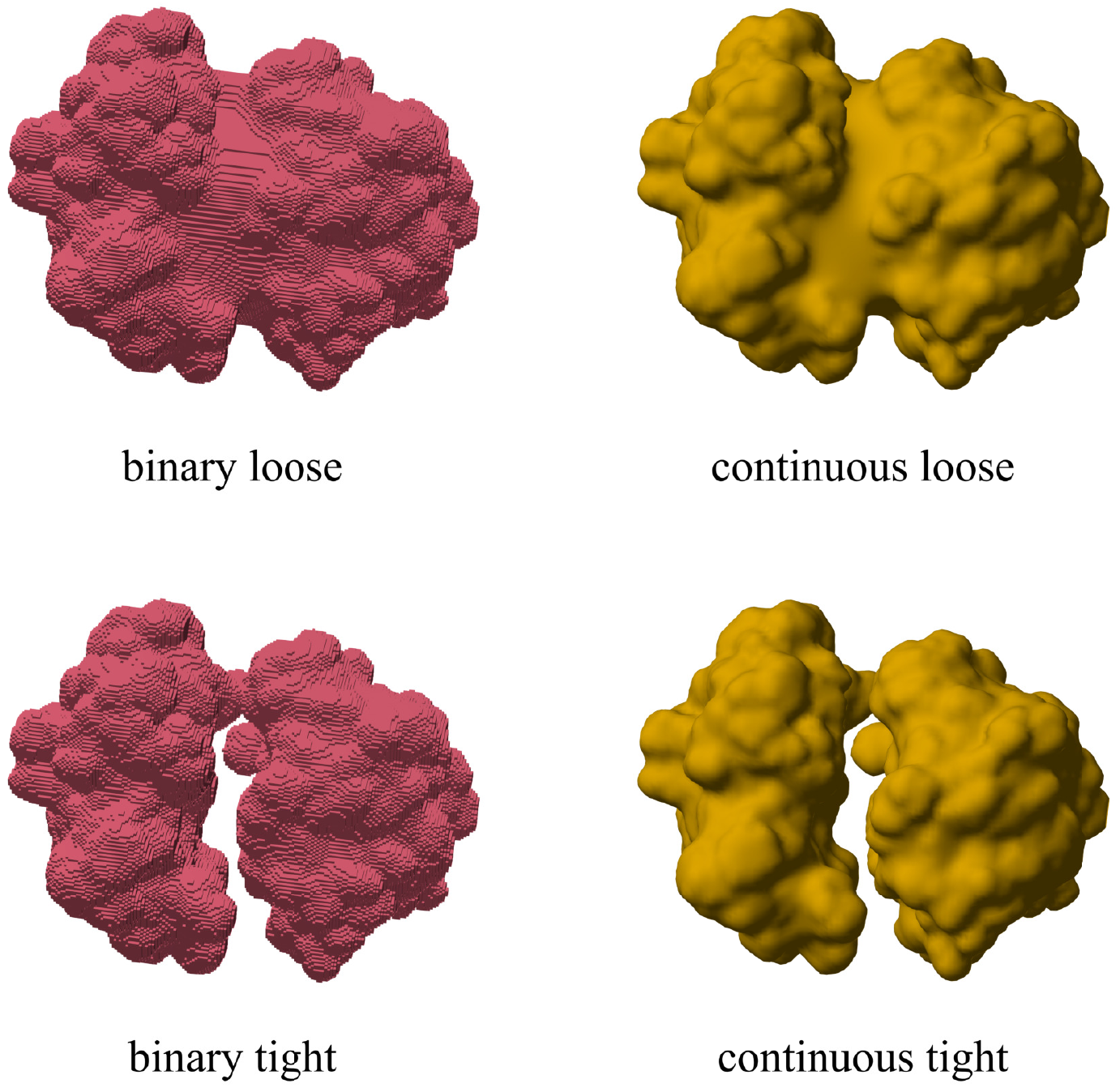}
\caption{Stable equilibrium solute−solvent interfaces of BphC obtained at d = 12 {\AA}. The computation box is $[-42.8425,42.8425]^3$ and $n=200$ for BLS-VISM, $n=106$ for CLS-VISM}
\label{f:bphc}
\end{figure}

In Fig.~\ref{f:p53mdm2} and ~\ref{f:bphc}, we compare the surfaces obtained from binary LSM and continuous LSM from tight or loose initial surfaces. In both BLS-VISM and CLS-VISM, the tight initial leads to a final state of two disjoint surfaces, which indicates that there is water between the two proteins, and this phenomenon is called wetting.
In contrast, the loose initial results in a final state of one connected interface, which indicateds that there is no water in between, and this phenomenon is called dewetting.
The ability of CLS-VISM to capture the dewetting effects of complex molecules is well-established \cite{wang_level-set_2012,zhou_variational_2014}. Here we demonstrate that the BLS-VISM preserves this characteristic feature of CLS-VISM.

In Table~\ref{t:protein}, we compare the energy and speed between BLS-VISM and CLS-VISM of the protein systems. 
For p53-MDM2 with either tight or loose fit, the relative error with respect to CLS-VISM is within 1\% for the surface energy, and within 5\% for the vdW and electrostatic energies. The relative error is 15\% for the total energy. 
For BphC with either tight or loose fit, the relative error with respect to CLS-VISM is within 3\% for the surface energy, and within 13\% for the vdW and electrostatic energies. The relative error is 50\% for the total energy. 
As mentioned before, the vdW and electrostatic energies are singular at the position of the atoms, and they are sensitive to the exact location of the interface. And the more atoms there are, the faster these two components go to positive or negative infinity as the surface get closer to the atoms.
The larger relative error in the total energy might be due to the cancellation of the large positive and negative energy components when computing the total energy, and this effect is more severe for larger system with more atoms. Nevertheless, due to the complex topology of the solvation interface and the large number of atoms, such discrepancies are expected and remain reasonably accurate when estimating solvation energy of protein systems \cite{zhang_coupling_2021}.

What's impressive is the speed of BLS-VISM, which is hundreds or thousands of times faster than CLS-VISM. The speedup is more significant as the number of atoms increases. This is because the work required in each flip of BLS does not depends on the number of atoms, but in CLS-VISM, when evolving the interface, the work required to evaluate the normal velocity \eqref{eq:vn} of the interface depends on the number of atoms. 
From the table, it seems that the speedup with loose initial is less impressive than that with tight fit. That's because in our original implementation of CLS-VISM with the loose initial, a coarse grid is used to speed up the evolution of the interface at the first stage. Similar idea can also be applied to BLS-VISM through a multi-resolution or adaptive grid, but we didn't pursuit along this line.

The significant speedup in computing VISM energy allow us to couple BLS-VISM with Monte Carlo (MC) Method to simulate protein binding. In \cite{zhang_coupling_2021}, the BLS-VISM is coupled with rigid-body MC to simulate binding of the p53-MDM2 system. During the MC simulation, the VISM energy needs to be evaluated millions of times, which is impossible with CLS-VISM. We note that in BLS-VISM with tight initial, the speed is no longer bottlenecked by the optimization process. Instead, the initialization process will take more time than the flipping process.

\begin{table}[!htpb]
\centering
\caption{Solvation free energy ($k_BT$) and computation time ($s$) for the protein systems}
\label{t:protein}
\medskip
\resizebox{\textwidth}{!}{
\begin{tabular}{|c|c|c|c|c|c|c|c|c|c|c|}
\hline\hline
\multirow{2}{*}{system}&\multicolumn{2}{|c|}{Surf} &\multicolumn{2}{c|}{vdW} &\multicolumn{2}{c|}{Elec} &\multicolumn{2}{c|}{total} &\multicolumn{2}{c|}{Time}\\
\cline{2-3} \cline{4-5}\cline{6-7} \cline{8-9} \cline{10-11}
&BLS & CLS & BLS & CLS& BLS& CLS& BLS& CLS& BLS &CLS\\
\hline
p53-MDM2 loose &  992.2 & 995.8 & -437.2 & -440.2 & -932.3 & -895.7 & -377.4 & -340.2 & 52.1 & 1510.1 \\
p53-MDM2 tight &  1019.4 & 1017.9 & -473.4 & -461.1 & -944.9 & -903.5 & -398.8 & -346.8 & 1.1 & 813.4 \\
BphC loose    &  2089.8 & 2032.0 & -1218.5 & -1083.0 & -1624.7 & -1456.6 & -753.3 & -507.6 & 52.2 & 25589.2 \\
BphC tight    &  2259.7 & 2215.4 & -1401.4 & -1240.6 & -1679.8 & -1528.2 & -821.4 & -553.3 & 1.2 & 14858.0 \\
\hline
\end{tabular}
}
\end{table}

\section{Conclusion}
\label{s:conclusion}
We have developed a fast binary level set method to minimize the VISM free-energy functional of solute-solvent interfaces. 
In the binary level set method, the interface is represented by a binary level set function that takes values $\pm 1$ on the solute or solvent region. 
A key component in our formulation is the approximation of surface area by convolution of indicator function with compactly supported kernel. 
The resulting discrete VISM energy is minimized by iteratively flipping the value of the binary level set function in a steepest descent fashion. 
As demonstrated by our numerical experiments, compared with the PDE-based level set method, the binary level set approach is hundreds of times faster, and still provide fairly accurate solvation energy. 
It also has the ability to capture different equilibrium solute-solvent interfaces. This is a characteristic feature of VISM and is important in biomolecular simulations.

Future works include estimating the curvature of an interface and incorporating the Poisson-Boltzmann theory of electrostatics into the binary level set framework.
Further performance gain might be achieved through an adaptive Cartesian grids.
We are also interested in applying our fast algorithm to simulate and study the binding and folding of biomolecules.

\appendix

\newcommand{\bn}{\mathbf{n}} 
\newcommand{\y}{\mathbf{y}} 
\newcommand{\hy}{\mathbf{\hat{y}}} 
\newcommand{\x}{\mathbf{x}} 
\newcommand{\z}{\mathbf{z}} 

\section{Integral Formulation of Surface Area}
\label{s:surfintegral}
In this section, we derive the integral formulation of surface area in the binary level set method \eqref{eq:area}. Let $\Gamma = \partial \Omega$ be a smooth hypersurface in $\R^d$.
Let $\phi$ be a signed distance function representing an interface with $\phi<0$ being the region $\Omega$ enclosed by the interface $\Gamma$,
\begin{equation}
   \phi(\x) = \begin{cases}
   -\inf_{\y\in\Omega}  |\x-\y| & \x \in \Omega \\
   \inf_{\y\in\Omega^c} |\x-\y| & \x \in \Omega^c \\
   \end{cases}.
\label{eq:signphi}
\end{equation}
Let $\Hv$ denote the one-dimensional Heaviside function, and $\delta$ the Dirac-delta function. We denote the $\theta$-level set of $\phi$ as $\Gamma_\theta=\{\x|\phi(\x) = \theta\}$ and $\Gamma_0 = \Gamma$. Let $\bn = \bn(\x)$ be the unit normal vector at $\x\in\Gamma$. As a property of the sign distance function, $\grad \phi(\x) = \bn$.

Then consider the following expression
\begin{equation}
\label{eq:V}
   \begin{aligned}
      V(\kr) & = \int_{\R^d} \Hv(-\phi(\x)) \int_{\R^d} \Hv(\phi(\y))K\left(\frac{|\x-\y|}{\kr}\right)d\y d\x\\
      & = \int_{\R^d} \Hv(-\phi(\x)) W(\x,\kr) d\x\\
      & = \int_{\{ -\kr \leq \phi(\x) \leq 0\}} W(\x,\kr) d\x\\
   \end{aligned}
\end{equation}
where 
\begin{equation}
\label{eq:W}
\begin{aligned}
   W(\x,\kr) &= \int_{\R^d} \Hv(\phi(\y))K\left(\frac{|\x-\y|}{\kr}\right)d\y\\
   & = \kr^d \int_{B_1(0)} \Hv(\phi(\kr \hy + \x))K(|\hy|)d\hy.
\end{aligned}
\end{equation}
The integrand in $V(\kr)$ is supported on $\{\x\mid -\kr \leq \phi(\x) \leq 0\}$ by the compactness of $K$.
Given $\z \in \Gamma_{-\theta}$, $\theta \in [0,\kr]$, there is a unique $\x \in \Gamma$ such that $\z = \x - \theta \bn$, as long as $\kr$ is smaller than the minimum radius of curvature.
Therefore we can reparametrize $\Gamma_{-\theta}$ by $\Gamma$.
By the coarea formula, we can write the volume \eqref{eq:V} as integral over level sets of $\phi$.
\begin{equation}
\begin{aligned}
      V(\kr) & = \int_{\{ -\kr \leq \phi(\x) \leq 0\}} W(\x,\kr) d\x\\
             & = \int_{0}^{\kr} \int_{\Gamma_{-\theta}}W(\x,\kr)d\x d\theta \\
             & = \int_{0}^{\kr} \int_{\Gamma} W(\z-\theta \bn , \kr) Y(\z,\theta)d\z d\theta \\
             & = \kr \int_{0}^{1} \int_{\Gamma} W(\z- \kr\theta \bn , \kr) Y(\z,\kr\theta)d\z d\theta \\
\end{aligned}
\end{equation}
Here, $Y(\z,s)$ is the Jacobian that accounts for the change of variables made in the reparametrization of $\Gamma_{-\theta}$ by $\Gamma$.
In two dimensions, $Y(\z,s) = 1 - H(\z)s$ and 
in three dimensions, $Y(\z,s) = 1 - H(\z)s +G(\z) s^2$. Here $H(\z)$ is the non-averaged mean curvature. $G(\z)$ is the Gaussian curvature. 
In $d$ dimensions, $Y(\z,s) = \prod_{i=1}^{d-1} (1-s \kappa_i(\z))$, where $\kappa_i$ is the $i$-th principal curvature of the hypersurface $\Gamma$ at $\z$ \cite{kublik_implicit_2013}.

\begin{figure}[!htpb]
\centering
\includegraphics[width=0.6\textwidth,keepaspectratio]{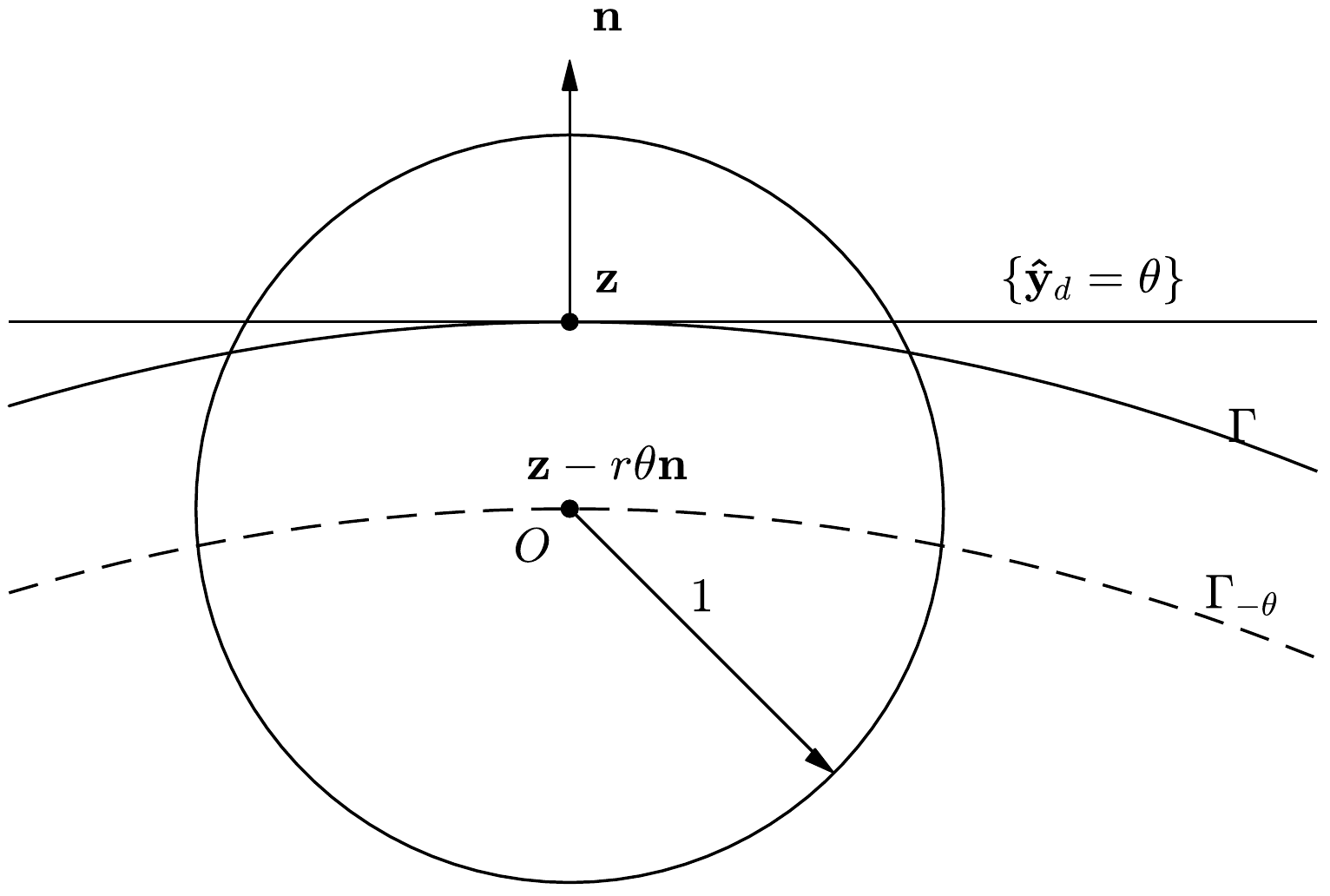}
\caption{Illustration of the term $W(\z-\kr\theta \bn,\kr)$.}
\label{f:w}
\end{figure}

We find the expansion of $W(\z-\kr\theta \bn,\kr)$ (see Fig. \ref{f:w}) with respect to $\kr$:
\begin{equation}
\begin{aligned}
   W(\z-\kr\theta \bn,\kr)& = \kr^d \int_{B_1(0)} \Hv(\phi(\z+\kr(\hy-\theta\bn)))K(|\hy|)d\hy\\
   &= \kr^d(a_0(\z,\theta) + \kr a_1(\z,\theta) + \bigO(\kr^2)).
\end{aligned}
\end{equation}
That leads to the expansion of $V(\kr)$. 
\begin{equation}
\begin{aligned}
V(\epsilon) &= \int_0^{1} \int_{\Gamma_0} W(\z-\epsilon \theta \bn,\epsilon)Y(\z,\epsilon \theta)\epsilon d\z d\theta\\
&= \kr^{d+1} \int_0^{1} \int_{\Gamma_0} a_0(\z,\theta) + \kr [a_1(\z,\theta) - H(\z)\theta a_0(\z,\theta)]d\z d\theta + \bigO(\kr^2)\\
&=\kr^{d+1} \left(v_0 + \kr v_1 + \bigO(\epsilon^2)\right)\\,
\end{aligned}
\end{equation}
where we define
\begin{equation}
   v_0 = \int_0^{1} \int_{\Gamma_0} a_0(\z,\theta) d\z d\theta,\quad
   v_1 = \int_0^{1} \int_{\Gamma_0} a_1(\z,\theta) - H(\z)\theta a_0(\z,\theta)d\z d\theta.
\end{equation}

Next, we simplifies the zero order term to find the constant $C_{K,\kr,d}$. By Taylor's expansion,
\begin{equation}
\begin{aligned}
      \phi(\z+\kr(\y-\theta \bn)) & = \phi(\z) + \grad \phi(\z) \cdot \kr (\y-\theta \bn) + \frac{1}{2} \kr(\y-\theta \bn) \cdot \grad^2\phi(\z) \kr(\y-\theta \bn) + \bigO(\kr^3) \\
      & = \kr \bn \cdot (\y-\theta\bn) + \frac{1}{2} \kr^2 (\y-\theta \bn) \cdot \grad^2\phi(\z) (\y-\theta \bn) + \bigO(\kr^3)\\
\end{aligned}
\end{equation}
Together with the fact that $\Hv$ is homogeneous of degree 0 (i.e. $\Hv(a\x) = \Hv(\x)$):
\begin{equation}
   \Hv(\phi(\z+\kr(\y-\theta \bn))) = \Hv( \bn \cdot (\y-\theta\bn) + \frac{1}{2} \kr (\y-\theta \bn) \cdot \grad^2\phi(\z) (\y-\theta \bn) + \bigO(\kr^2)),
\end{equation}
we have 
\begin{equation}
   \lim_{\kr\to0} \Hv(\phi(\z+\kr(\y-\theta \bn))) = \Hv(\bn\cdot(\y-\theta \bn)),
\end{equation}
where $\{\y \mid \bn \cdot (\y-\theta \bn) = 0\} $ is the hyperplane passing through $\theta \bn$ with normal vector $\bn$. Therefore $\{\y \mid \bn \cdot (\y-\theta \bn) \geq 0\} $ is the side of the hyperplane in the direction of $\bn$. Without loss of generality, we can assume $\bn = e_d$.
\begin{equation}
\begin{aligned}
   a_0(\z,\theta) & =  \lim_{\kr\to 0} W(\z-\kr\theta\bn,\kr)\\
   & = \int_{ B_1(0)} \Hv(\bn\cdot(\hy-\theta \bn))K(|\hy|)d\hy\\
   & = \int_{ B_1(0) \cap \{ \bn \cdot (\hy-\theta \bn) \geq 0\}}  K(|\hy|)d\hy\\
   & = \int_{ B_1(0) \cap \{\hy_d \geq \theta\} }K(|\hy|)d\hy\\
\end{aligned}
\end{equation}
$\hy_d$ is the d-th coordinate of $\hy$.
Note $a_0(\z,\theta) = a_0(\theta)$ does not depend on $\z$. Hence, 
\begin{equation}
\begin{aligned}
   v_0 & = \int_0^1\int_{\Gamma_0}a_0(\theta)dAd\theta\\
   & = \left(\int_0^1 a_0(\theta)d\theta\right) {\rm Area} (\Gamma)\\
\end{aligned}
\end{equation}

The constant can be further simplified by writing the integral in polar coordinates $r = |\y| \in (0,1] $ and $\x = \y/|\y| \in S^{d-1}$
\begin{equation}
\begin{aligned}
   \int_0^1 a_0(\theta)d\theta &=  \int_0^1 \int_{ B_1(0) } \Ind_{\{\y_d \geq \theta\}}  K(|\y|)d\y d\theta \\
   & = \int_0^1 \int_0^1 \int_{ S^{d-1}} \Ind_{\{r\x_d \geq \theta\}}  K(r)r^{d-1} d\x dr d\theta \\
   & = \int_0^1 \int_{ S^{d-1}} \left(\int_0^1 \Ind_{\{r\x_d \geq \theta\}}  d\theta\right) K(r)r^{d-1} d\x dr \\
   & = \int_0^1 \int_{ S^{d-1}} \Ind_{\{\x_d>0\}} r\x_d K(r)r^{d-1} d\x dr \\
   & = \left(\int_{ S^{d-1}\cap\{\x_d>0\}} \x_d d\x \right) \int_0^1 K(r)r^d dr \\
   & = C_d \int_0^1 K(r)r^d dr \\
\end{aligned}   
\end{equation}
where
\begin{equation}
   C_d = \frac{2 \pi^{\frac{d-1}{2}}} {(d-1)\Gamma(\frac{d-1}{2})}
\end{equation}
For $d = 3$, $C_d = \pi$. For $d=2$, $C_d = 2$.



For the first order term, we use the identity that $0 = \grad (\bn\cdot\bn) = 2\grad^2\phi(\z)\bn$. And without loss of generality, we can assume $\bn = e_d$, and hence $[\grad^2\phi(\z)]_{dd}=0$

\begin{equation}
   \begin{aligned}
   a_1(\z,\theta) & = \frac{d}{d\kr}\Bigr|_{\kr=0} \int_{B_1(0)} \Hv(\phi(\z+\kr(\y-\theta\bn)))K(|\y|)d\y \\
   & =\int_{B_1(0)}  \frac{d}{d\kr}\Bigr|_{\kr=0} \Hv(\phi(\z+\kr(\y-\theta\bn)))K(|\y|)d\y \\
   & = \int_{B_1(0)} \delta(\bn \cdot (\y-\theta\bn)) \frac{1}{2}(\y-\theta\bn)\cdot \grad^2\phi(\z) (\y-\theta\bn)  K(|\y|)d\y \\
   & =  \frac{1}{2} \int_{B_1(0) \cap \{\y_d=\theta\}}(\y-\theta\bn)\cdot \grad^2\phi(\z) (\y-\theta\bn)  K(|\y|)d\y \\
   & =  \frac{1}{2} \int_{B_1(0) \cap \{\y_d=\theta\}}\y\cdot \grad^2\phi(\z) \y  K(|\y|)d\y \\
   & =  \frac{1}{2} \int_{B_1(0) \cap \{\y_d=\theta\}}\y\cdot \grad^2\phi(\z) \y  K(|\y|)d\y \\
   & = \frac{1}{2}  \sum_{i,j=1}^d \int_{B_1(0) \cap \{\y_d=\theta\}} \y_i\y_j [\grad^2\phi(\z)]_{ij} K(|\y|)d\y \\
   \end{aligned}
\end{equation}
The region of integration is a spherical section of the unit ball. Suppose $i\neq j$, and without loss of generality, suppose $j\neq d$, then the integrand is an odd function with respect to $\y_j$. By symmetry of the region, the integral is 0. Therefore

\begin{equation}
   \begin{aligned}
   a_1(\z,\theta) & = \frac{1}{2}  \sum_{i,j=1}^d \int_{B_1(0) \cap \{\y_d=\theta\}} \y_i\y_j [\grad^2\phi(\z)]_{ij} K(|\y|)d\y \\
   & = \frac{1}{2} \left( \sum_{i=1}^d [\grad^2\phi(\z)]_{ii} \right) \int_{B_1(0) \cap \{\y_d=\theta\}} \y_1^2 K(|\y|)d\y \\
   & = \frac{1}{2} H(z) \int_{B_1(0) \cap \{\y_d=\theta\}} \y_1^2 K(|\y|)d\y \\
   & = H(z) \hat{a}_1(\theta) \\
   \end{aligned}
\end{equation}
Note that
\begin{equation}
   \begin{aligned}
   \int_0^1\theta a_0(\theta)d\theta &=\int_0^1 \int_{B_1(0) \cap \{\y_d\geq\theta\} } \theta K(|\y|) d\y d\theta\\
   & =  \int_{B_1(0) \cap \{\y_d\geq 0\} } \int_0^1 \Hv(\y_d-\theta)\theta K(|\y|) d\theta d\y\\
   & =  \int_{B_1(0) \cap \{\y_d\geq 0\} } \left(\int_0^{\y_d} \theta d\theta \right) K(|\y|) d\y\\
   & =   \frac{1}{2}  \int_{B_1(0) \cap \{\y_d\geq 0\} }{\y_d}^2 K(|\y|) d\y\\
   & =  \int_0^1 \hat{a}_1(\theta)d\theta\\
   \end{aligned}
\end{equation}

Therefore
\begin{equation}
   v_1 = \int_{\Gamma_0} \int_0^{1}  a_1(\z,\theta) - \theta H(\z) a_0(\z,\theta)d\theta dA = 0
\end{equation}
and 
\begin{equation}
\label{eq:final area}
   {\rm Area}(\Gamma) = C_{K,\kr,d} V(\kr) + \bigO(\kr^2)
\end{equation}
where 
\begin{equation}
   C_{K,\kr,d} = \left(\kr^{d+1} C_d\int_0^1 K(\kr)\kr^dd\kr\right)^{-1}   
\end{equation}

\section{Numerical Approximation of the Surface Area}
\label{s:surfnumeric}
\newcommand{\rall}{\Omega} 
\newcommand{\rin}{\Omega_{in}}
\newcommand{\rout}{\Omega_{out}}
\newcommand{\ha}{\hat{\alpha}}
\newcommand{\hb}{\hat{\beta}}
\renewcommand{\a}{\alpha}
\renewcommand{\b}{\beta}

In this section, we analyze the numerical error of computing the integral \eqref{eq:V} with midpoint rule. 
Let $\Gamma = \partial \rin$ be a smooth compact hypersurface contained in a cube $\rall \subset \R^d$. $\rall = \rin \cup \rout$.
Beside assuming $K(\x)$ to be a radially symmetric compact kernel with unit radius, we also assuming that $K$ is twice continuously differentiable, so that $K$ goes to 0 smoothly at $|\x|=1$.

Suppose $\rall$ is covered with a uniform Cartesian grid of size $h$, with $n$ grid cells in each dimension. 
Let $c_i$ be the grid cell centered at $\x_i$ with side length $h$, where $i$  is a multi-index $(i_1,\dots,i_d)$, with $i_k \in \{1,\dots,n\}$, $1\leq k \leq n$.
With a slight abuse of notation, we also define the index set $\rin$: $i\in \rin$ if $\x_i \in \rin$. Similarly for $\rall$ and $\rout$. We also define the index set $I = \{i \mid c_i \cap \Gamma \neq \emptyset\}$. $I$ contains all the ``interface points'', whose grid cells touch the interface. $I^c$ contains all the ``interior points'', whose grid cells are completely in $\rin$ or $\rout$. 

Recall the midpoint rule for numerical integration. Suppose $f$ is twice continuously differentiable. Let $\x^*$ be the center of a $d$ dimensional cube of side length $h$. Then the approximation on the cube has error $\bigO(h^{d+2})$.
\begin{equation}
\label{eq:midfull}
   \begin{aligned}
      \int_{[0,h]^d} f(\x) d\x & = \int_{[0,h]^d} f(\x^*) + (\x-\x^*)\cdot\grad f(\x^*) + \bigO(h^2) d\x\\
      &=h^d f(\x^*) + \bigO(h^{d+2})\\
   \end{aligned}
\end{equation}

If the region is not a cube but some region $\omega \subset [0,h]^d$, and the volume of the region is known, then our error is $\bigO(h^{d+1})$
\begin{equation}
\label{eq:midexact}
   \begin{aligned}
      \int_{\omega} f(\x) d\x & = \int_{\omega} f(\x^*) + \bigO(h) d\x \\
      &={\rm vol}(\omega) f(\x^*) + \bigO(h^{d+1})\\
   \end{aligned}
\end{equation}

However, if the exact volume of the region $\omega$ is unknown, then the following approximation has error $\bigO(h^{d})$
\begin{equation}
\label{eq:midbinary}
   \begin{aligned}
      \int_{\omega} f(\x) d\x & = h^d \Ind_\omega (\x^*) f(\x^*) + \bigO(h^d) \\
   \end{aligned}
\end{equation}


Recall the formula of $W(\x,\kr)$ \eqref{eq:W} and $V(\kr)$ \eqref{eq:V}:
\begin{equation}
   \begin{aligned}
      W(\x,\kr) &= \int_{\rout} K\left(\frac{|\x-\y|}{\kr}\right)d\y\\
      V(\kr) & = \int_{\rin} W(\x,\kr)  d\x\\
   \end{aligned}
\end{equation}

Let $\ha_i = {\rm vol}(c_i \cap \rin)$,  $\hb_i = {\rm vol}(c_i \cap \rout)$. These are the exact volume fractions of the grid cells. We define the intermediate quantity $\hat{V}_h(\kr)$, which is composite midpoint rule approximation of $V(\kr)$ using exact volume fraction. 
\begin{equation}
   \hat{V}_h(\kr) = \sum_{i,j} h^{2d} \ha_i  \hb_i K_{r,ij},
\end{equation}
where the summation is taken over all ordered tuples of indices $(i,j)\in \rall \times \rall$, and
\begin{equation}
   K_{r,ij} = K\left(\frac{|\x_i-\x_j|}{\kr}\right),
\end{equation}

To look at the error of $\hat{V}_h(\kr)$, firstly, we approximate $W(\x_i,\kr)$ by composite midpoint rule with exact volume fraction. 
Consider $j\in \rall$, if $j$ is an interface point, the error is $\bigO(h^{d+1})$ by \eqref{eq:midexact},
and the number of interface points in a ball of radius $\kr$ is $\bigO(\kr^{d-1}h^{1-d})$. Hence the total error due to interface points is $\bigO(r^{d-1} h^2)$.
If $j$ is an interior point, we have two cases:
if $\x_j-\x_i>\kr$, then the error is 0. 
if $\x_j-\x_i\leq\kr$, then the error is $\bigO(h^{d+2})$ by \eqref{eq:midfull}, and the numer of interior points in a ball of radius $\kr$ is $\bigO(\kr^{d}h^{1-d})$.
Here we make use of the assumption that $K$ goes to 0 smoothly at the edge of the kernel.
The total error due to interior points is $\bigO(r^d h^2)$.
Therefore, in total we have
\begin{equation}
    W(\x_i,\kr) = \sum_{j} h^d \hb_{j} K_{r,ij} + \bigO(\kr^{d-1}h^2).
\end{equation}

Next, we approximate $V(\kr)$ by composite midpoint rule with exact volume fraction. The number of points near the interface is given by $|I| = \bigO(h^{1-d})$. Since $W(\x_i,\kr) = \bigO(\kr^d)$, the error on individual interface point is $\bigO(\kr^d h^{d+1})$. The total error is given by

\begin{equation}
   \begin{aligned}
   V(\kr) & = \sum_{i} h^d \ha_i W(\x_i,\kr) + \bigO(h^2\kr^{d})\\
   \end{aligned}
\end{equation}

Plug in the estimation of $ W(\x_i,\kr)$, and recall that $W(\x_i,\kr)$ is nonzero for points in a tubular neighborhood of the interface. Let $T = \{ i\mid -\kr \leq \phi(\x_i) \leq 0\}$ and $|T| = \bigO(\kr/h^d)$. Hence
\begin{equation}
   \begin{aligned}
   V(\kr) &= \sum_{i} \sum_{j} h^{2d} \ha_i  \hb_i K_{r,ij} + \bigO( h^d (\kr/h^d) (h^2\kr^{d-1}) + h^2\kr^{d})\\
   & = \hat{V}_h(\kr)+ \bigO( h^2\kr^{d})\\
   \end{aligned}
\end{equation}


We also define $V_h(\kr)$, which is composite midpoint rule approximation of $V(\kr)$ without knowing the exact volume fraction:
\begin{equation}
   V_h(\kr) = \sum_{i,j} h^{2d} \a_i \b_{j} K_{r,ij}
\end{equation}
where $\a_i = 1$ if $i \in \rin$, and $\a_i = 0$ if $i \in \rout$. $b_i = 1 - a_i$.

\begin{figure}[!htpb]
\centering
\includegraphics[width=0.6\textwidth,keepaspectratio]{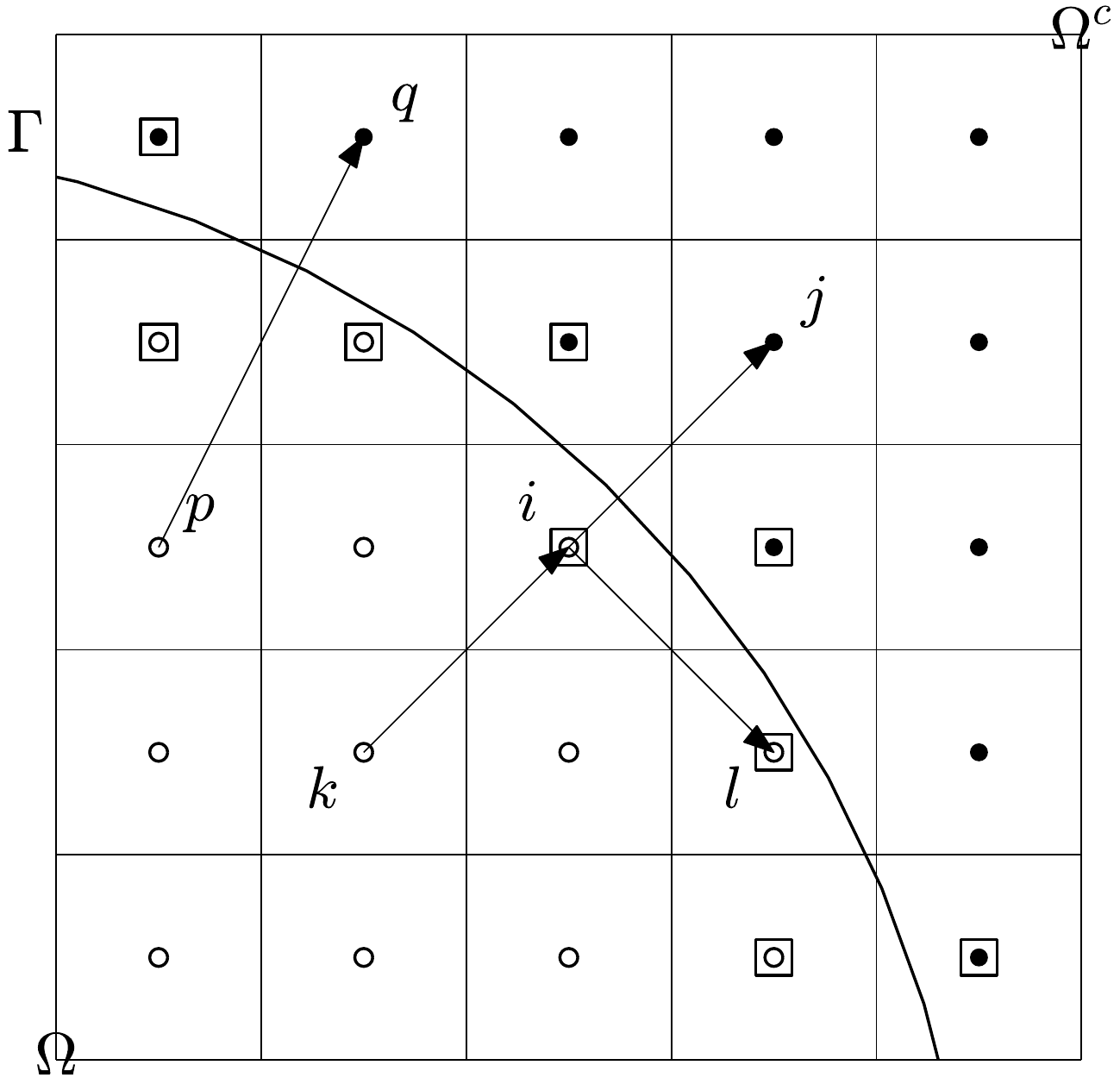}
\caption{Illustration of different cases. Squared grid points are in the set $I$, their grid cell intersects with the interface $\Gamma$. $(p,q) \in Z_1$. $(i,j)$ and $(k,i)$ are paired and belong to $Z_2$. $(i,l) \in Z_3$}
\label{f:pair}
\end{figure}

We analyze the difference between $V_h(\kr)$ and $\hat{V}_h(\kr)$
\begin{equation}
   V_h(\kr) - \hat{V}_h(\kr) = \sum_{i,j} h^{2d} (\a_i \b_{j} - \ha_i \hb_{j}) K_{r,ij}\\
\end{equation}
by partitioning $\rall \times \rall$ into three disjoint subsets $Z_1$,$Z_2$ and $Z_3$, cf Fig. \ref{f:pair}.

\textbf{Case 1}: $Z_1 = I^c \times I^c$.
Both $i$ and $j$ are interior points. Then $\ha_i = 1$ if $i\in \rin$, and 0 otherwise. So $\ha_i = \a_i$ and $\hb_i = \b_i$. On this set, the difference between the two methods is 0.

\textbf{Case 2}: 
We define $R_i(j) = 2 i-j$, which is the point reflection of $j$ with respect to $i$. Let $\hat{\phi}_i = {\rm sgn}(\phi(\x_i))$ be the sign of the level set function at $\x_i$. So the condition $\hat{\phi}_{i} \hat{\phi}_{j} = -1$ means that $\x_i$ and $\x_j$ are on opposite sides of the interface. 
And define the following two sets that are disjoint,
\begin{equation}
\begin{aligned}
S_1 &= \{(i,j)\mid i\in I, j\in I^c, R_i(j) \in I^c, \hat{\phi}_{R_i{(j)}} \hat{\phi}_{j} = -1\},\\
S_2 &= \{(i,j)\mid i\in I^c, j\in I, R_j(i) \in I^c, \hat{\phi}_{R_j{(i)}} \hat{\phi}_{i} = -1\}.\\
\end{aligned}
\end{equation}
In words, for tuples $(i,j)$ in $S_1$, $i$ is an interface point, $j$ is an interior point, the reflection of $j$ is also an interior point on the opposition side of the interface. By definition, the map $(i,j)\mapsto(R_i(j),i)$ is a bijection between $S_1$ and $S_2$. For case 2, $Z_2 = S_1 \cup S_2$.

Given that $(i,j)\in Z_2$, we can pair up the error, 
\begin{equation}
   \begin{aligned}
   & \sum_{(i,j)\in S_1\cup S_2} h^{2d} (\a_i \b_{j} - \ha_i \hb_{j}) K_{r,ij} \\
   = &\sum_{(i,j)\in S_1} h^{2d} \left[(\a_i \b_{j} - \ha_i \hb_{j}) + (\a_{k} \b_{i} - \ha_{k} \hb_{i})\right] K_{r,ij}\\
   \end{aligned}
\end{equation}
where $k = R_i(j)$. Because both $j$ and $k$ are interior points, $\hb_j = \b_j$ and $\ha_k = \a_k$. Because $j$ and $k$ are in different side, $\a_k = 1 - \a_j = \b_j$. Together with the fact that $\a_i + \b_i = 1$ and $\ha_j + \hb_j = 1$, we have
\begin{equation}
   \begin{aligned}
   (\a_i \b_{j} - \ha_i \hb_{j}) + (\a_{k} \b_{i} - \ha_{k} \hb_{i}) = \a_i \b_j - \ha_i \b_j + \b_j\b_i-\b_j \hb_i = 0.
   \end{aligned}
\end{equation}
Hence, on case 2, the difference between the two methods also is 0.

\textbf{Case 3}: 
$Z_3$ contains all the remaining tuples, where there is no cancellation that can be exploited, so $\a_i \b_{j} - \ha_i \hb_{j} = \bigO(1)$. Our goal is to bound the size of $Z_3$. 
Given $i\in I$, we define the following index set:
\begin{equation}
   Q_{r,i} = \{|\x_i - \x_j| \leq \kr\} \cap \left( \{j\in I\}\cup \{j\in I^c, R_i(j)\in I \text{ or } \hat{\phi}_{R_i{(j)}} \hat{\phi}_{j} = 1 \} \right)
\end{equation}
In words, given that $i$ is an interface point, among all the points $\x_j$ that are within kernel radius of $\x_i$, either $j$ is also an interface point, or $j$ is an interior point whose reflection is an interface point or an interior point on the same side as $i$. Then $Z_3$ can be written as
\begin{equation}
   \begin{aligned}
      Z_3 & = \{(i,j) \mid i\in I, j\in Q_{r,i} \} \cup \{(i,j) \mid j\in I, i\in Q_{r,j}\} \\
   \end{aligned}
\end{equation}

\begin{figure}[!htpb]
\centering
\includegraphics[width=0.6\textwidth,keepaspectratio]{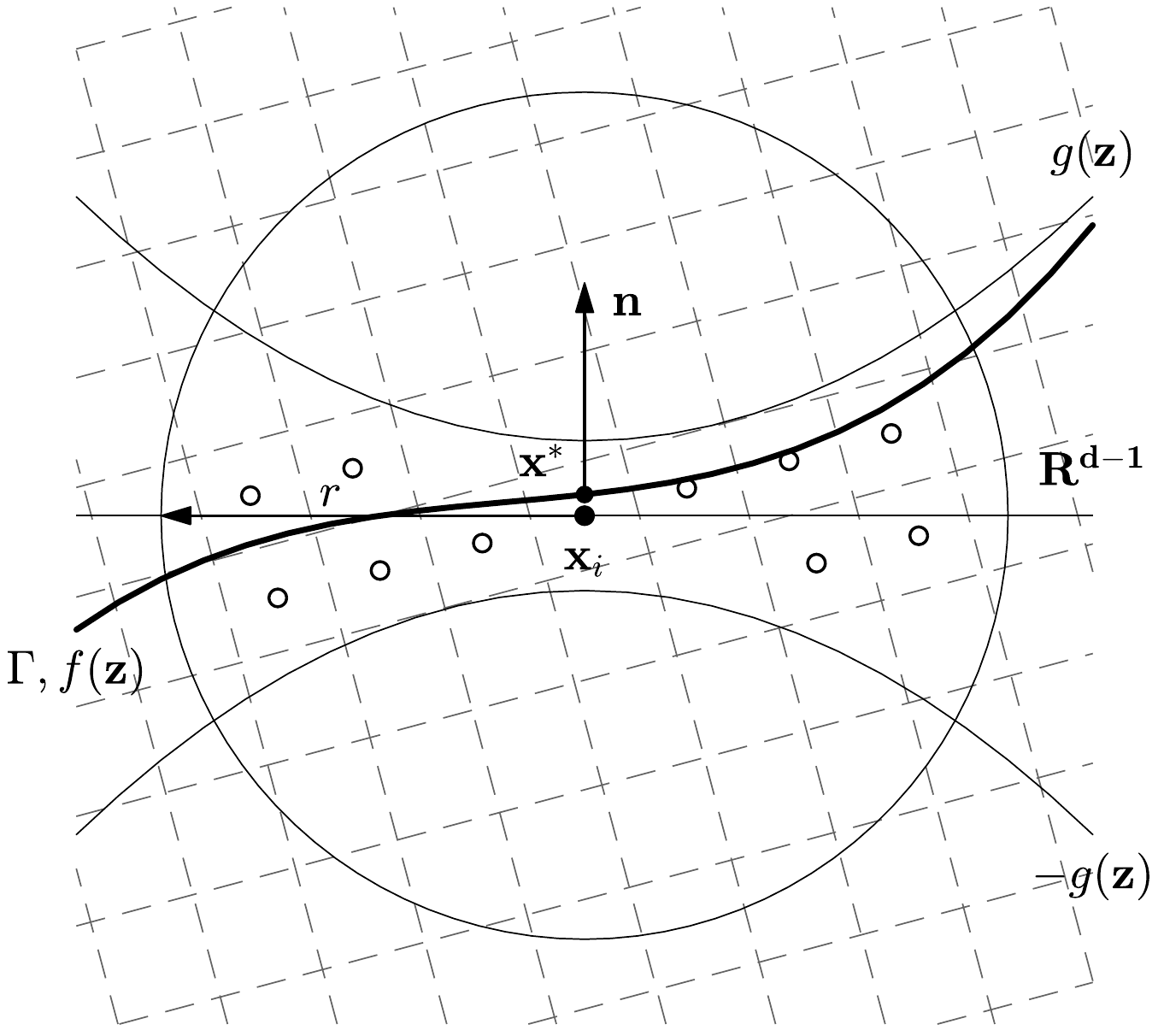}
\caption{Estimation of $|Q_{r,i}|$. $\x_i$ is an interface point. $\x^*$ is the closest point to $\x_i$ on $\Gamma$. In local coordinate, $\Gamma$ is the graph of the function $f(\z)$. Points in $Q_{r,i}$ are circled and bounded between the quadratic functions $g(z)$ and $-g(z)$}
\label{f:case3}
\end{figure}

We use the fact that an embedded hypersurface can be locally approximated by the graph of a quadratic function. Given $\x_i\in I$, we can find the closest point $\x^*$ on the interface $\Gamma$, so that $\x_i = \x^* - s \bn(\x^*)$ for some $s$, and $|s|\leq \sqrt{d} h/2$, which is half the length of the diagonal of a grid cell. We can create a local coordinate where $\x_i$ is the origin, and the first $d-1$ coordinate is parallel to the tangent plane of $\Gamma$ at $\x^*$. Then locally, $\Gamma$ can be parametrized as a graph $(\z,f(\z))$, $\z\in\R^{d-1}$, 
\begin{equation}
   f(\z) = s + \frac{1}{2} \left(\kappa_1 z_1^2+\cdots+\kappa_{d-1} z_{d-1}^2 \right) + \bigO(|\z|^3), 
\end{equation}
where $\kappa_1,\dots,\kappa_{d-1}$ are the principal curvatures of $\Gamma$ at $\x^*$. 
So $f(\z)$ will be bounded above by some quadratic function $g(\z)$,
\begin{equation}
   g(\z) = |s|+ \bigO(|\z|^2), 
\end{equation}
And $-g(\z)<f(\z)<g(\z)$. 
Hence, if $j \in Q_{r,i}$, $\x_i$ and $\x_j$ will both be in the region between $g(\z)$ and $-g(\z)$. Otherwise $(i,j) \in Z_2$.

By the compactness of the kernel, we only need to consider $\z\in B(\kr)$, a $d-1$ dimensional ball of radius $\kr$.
The volume of the region under the graph is bounded by
\begin{equation}
   \int_{B(r)} |s|+ \bigO(\z^2)d\z = \bigO(h\kr^{d-1} + \kr^{d+1}).
\end{equation}
Hence 
\begin{equation}
   |Q_{r,i}| =  \bigO \left(h^{1-d}\kr^{d-1} + h^{-d}\kr^{d+1}\right)
\end{equation}
and the size of $Z_3$ can be bounded
\begin{equation}
   |Z_3| \leq 2 |I| |Q_{r,i}| = \bigO\left( \kr^{d-1}h^{2(1-d)} + \kr^{d+1}h^{1-2d} \right).
\end{equation}
And we obtained the error between $\hat{V}_h(\kr)$ and $V_h(\kr)$:
\begin{equation}
\begin{aligned}
   \hat{V}_h(\kr) & = V_h(\kr) + \bigO(|Z_3|) h^{2d}\\
   & = V_h(\kr) + \bigO(h^2 r^{d-1} + hr^{d+1})\\
\end{aligned}
\end{equation}
Hence the error between $\hat{V}_h(\kr)$ and $V(\kr)$ is given by:
\begin{equation}
\begin{aligned}
   V(r) & = \hat{V}_h(\kr) + \bigO(h^2\kr^d) \\
   & = V_h(\kr) + \bigO(h^2 r^{d-1} + hr^{d+1} + h^2\kr^d)
\end{aligned}
\end{equation}
Recall that $C_{K,r,d} = \bigO(\kr^{-(d+1)})$, we have 
\begin{equation}
   \begin{aligned}
   {\rm Area}(\Gamma) &= C_{K,\kr,d} V(\kr) + \bigO(\kr^2)\\
   & = C_{K,\kr,d} V_h(\kr) +  \bigO(h^2/r^2 + h + h^2/\kr + \kr^2)
   \end{aligned}
\end{equation}
Therefore, for $\kr \sim \sqrt{h}$, 
we obtained a first order approximation of the surface area.
\begin{equation}
      {\rm Area}(\Gamma) = C_{K,\kr,d} V_h(\kr) + \bigO(h)
\end{equation}

\section{Integration in outside box region}
\label{s:outside}
For details on how to compute the integral $G_{vdW}[\Gamma]$ and $G_{elec}[\Gamma]$ outside of the computational box, see \cite{cheng_level-set_2010}. The idea is to partition the region outside of the computational box $\Omega^c$,  write the integral in different partition in cylindrical or spherical coordinate, integrate analytically in two of the dimensions, and finally compute the one-dimensional integral with composite Simpson's rule.

Here we discuss another method that make use of parallel computing. Let $\Omega$ be a cube of equal side in $\mathbb{R}^3$, centered at the origin, and $S$ be the sphere inscribed in $\Omega$. Then we can write the integral in spherical coordinate
\begin{equation}
\begin{split}
   \int_{\Omega_c} f(x)dx & = \int_{S^c} \Ind_{\Omega}(x,y,z) f(x,y,z) dxdydz \\
   & = \int_{0}^{\pi} \int_{0}^{2\pi} \int_{R}^{\infty} \hat{\Ind}_{\Omega}(\theta,\varphi,r) \hat{f}(\theta,\varphi,r) r^2 \sin\theta dr d\theta d\varphi\\
   & = \int_{0}^{\pi} \int_{0}^{2\pi} \int_{0}^{1/R} \hat{\Ind}_{\Omega}(\theta,\varphi,\frac{1}{\rho}) \hat{f}(\theta,\varphi,\frac{1}{\rho}) \rho^{-4} \sin\theta d\rho d\theta d\varphi\\
\end{split}
\end{equation}
In the second equality, we change from Cartesian coordinate to Polar coordinate. In the last equality, we perform a change of variable $\rho=r^{-1}$.

We apply the midpoint rule
\begin{equation}
\int_{\Omega_c} f(x)dx \approx \sum_{i,j,k} \hat{\Ind}_{\Omega}(\theta_j,\varphi_k,\frac{1}{\rho_i}) \hat{f}(\theta_j,\varphi_k,\frac{1}{\rho_i}) \rho^{-4} h^3 + \bigO(h)
\end{equation}

Note that the summand can be computed in parallel. Through experiment we found that our implementation of this method (using OpenCL and OpenMP) achieves desirable accuracy in shorter time compared with the method in \cite{cheng_level-set_2010}.

\bibliographystyle{unsrt}
\bibliography{binaryls}

\end{document}